\documentclass[11pt, a4paper, UKenglish]{article}

\usepackage[utf8]{inputenc}
\usepackage[english]{babel}

\usepackage{amsmath}
\usepackage{amssymb}
\usepackage{amsthm}
\usepackage{tikz}
\usepackage{tikz-cd}
\usepackage{todonotes}
\usepackage{hyperref}
\usepackage{footmisc}

\usepackage{bbm}
\usepackage{afterpage}
\usepackage{mathrsfs}

\newcommand{\M}{\mathcal{M}}
\newcommand{\Z}{\mathbb{Z}}
\newcommand{\Q}{\mathbb{Q}}
\newcommand{\R}{\mathbb{R}}

\newcommand{\KK}{\mathbf{K}}

\newcommand{\al}{\alpha}
\newcommand{\be}{\beta}

\newcommand{\ga}{\gamma}

\renewcommand{\O}{\mathcal{O}}

\renewcommand{\P}{\mathbb{P}}

\newcommand{\tor}{\mathrm{Tor}}

\DeclareMathOperator{\colim}{colim}

\newcommand{\nat}{\mathrm{Nat}}
\renewcommand{\hom}{\mathrm{Hom}}

\renewcommand{\mod}{\mathrm{Mod}}
\newcommand{\lmod}{\mathrm{LMod}}
\newcommand{\rmod}{{\mathrm{RMod}}}

\newcommand{\opend}{{\mathscr{E}\mathrm{nd}}}

\renewcommand{\top}{\mathrm{Top}}
\newcommand{\spec}{\mathrm{Sp}}
\newcommand{\topop}{\mathrm{TopOp}}

\newcommand{\ab}{\mathrm{Ab}}

\newcommand{\gf}{\mathcal{GF}}

\newcommand{\C}{\mathscr{C}}

\newcommand{\id}{\mathrm{id}}
\newcommand{\im}{\mathrm{im}}

\newcommand{\Ab}{\mathcal{A}\mathrm{b}}

\newcommand{\ho}{\mathrm{Ho}}
\newcommand{\Ho}{{\mathrm{Ho}}}

\DeclareMathOperator{\hocolim}{hocolim}

\newcommand{\calg}{\mathrm{CAlg}}

\newcommand{\smp}{\otimes}

\newcommand{\smpl}{\otimes^\mathbb{L}}
\newcommand{\smplr}{\underset{R}{\otimes}^\mathbb{L}}
\newcommand{\sph}{\mathbb{S}}

\renewcommand{\vee}{\oplus}
\renewcommand{\bigvee}{\bigoplus}

\newcommand{\Map}{\mathrm{Map}}

\newcommand{\mU}{\mathrm{MU}}
\newcommand{\MU}{\mathbf{MU}}
\newcommand{\BP}{\mathbf{BP}}

\newcommand{\1}{\mathbbm{1}}

\newcommand{\epi}{\underline{\pi}}
\newcommand{\res}{\mathrm{res}}
\newcommand{\tr}{\mathrm{tr}}

\newcommand{\burn}{\mathbf{A}}
\newcommand{\gl}{\mathrm{gl}}

\renewcommand{\gf}{\mathcal{GF}}

\newcommand{\U}{\mathcal{U}}

\newcommand{\F}{{\mathscr{F}}}

\newcommand{\E}{{\mathbb{E}}}
\newcommand{\G}{{\mathbb{G}}}

\renewcommand{\int}{\mathrm{int}}

\theoremstyle{plain}
\newtheorem{theorem}[equation]{Theorem}

\newtheorem*{theorema}{Theorem A}
\newtheorem*{theoremb}{Theorem B}

\newtheorem{prop}[equation]{Proposition}
\newtheorem{lemma}[equation]{Lemma}
\newtheorem{cor}[equation]{Corollary}

\theoremstyle{definition}
\newtheorem{mydef}[equation]{Definition}
\newtheorem{example}[equation]{Example}
\newtheorem{remark}[equation]{Remark}

\newtheorem{construction}[equation]{Construction}

\numberwithin{equation}{section}

\begin{document}

\title{Realising $\pi_\ast^eR$-algebras by global ring spectra}


\author{J.M. Davies\footnote{\url{j.m.davies@uu.nl}}}

\maketitle

\begin{abstract}
We approach a problem of realising algebraic objects in a certain universal equivariant stable homotopy theory; the global homotopy theory of Schwede \cite{s}. Specifically, for a global ring spectrum $R$, we consider which classes of ring homomorphisms $\eta_\ast\colon\pi_\ast^e R\rightarrow S_\ast$ can be realised by a map $\eta\colon R\rightarrow S$ in the category of global $R$-modules, and what multiplicative structures can be placed on $S$. If $\eta_\ast$ witnesses $S_\ast$ as a projective $\pi_\ast^e R$-module, then such an $\eta$ exists as a map between homotopy commutative global $R$-algebras. If $\eta_\ast$ is in addition \'{e}tale or $S_0$ is a $\Q$-algebra, then $\eta$ can be upgraded to a map of $\E_\infty$-global $R$-algebras or a map of $\G_\infty$-$R$-algebras, respectively. Various global spectra and $\E_\infty$-global ring spectra are then obtained from classical homotopy theoretic and algebraic constructions, with a controllable global homotopy type.
\end{abstract}

\tableofcontents

\newpage

\section*{Introduction}
\addcontentsline{toc}{section}{Introduction}
A key feature of the stable homotopy category is the interplay between algebra and homotopy theory. In this article we explore a variant of the following realisation problem:

\begin{center}\emph{Given a ring spectrum $R$, when does a map of graded rings $\pi_\ast R\rightarrow S_\ast$ come from a map of $R$-module spectra $R\rightarrow S$?}\end{center}

\begin{example}\label{emexample}
If $R$ is an Eilenberg--Mac~Lane spectrum and $S_\ast$ is concentrated in degree zero, then the answer is ``always''.  One way to see this is to recognise the full subcategory of the ($\infty$-) category of $R$-module spectra spanned by Eilenberg--Mac~Lane spectra as the (nerve of the) category of $\pi_0 R$-modules; see \cite[Proposition 7.1.1.13(3)]{ha}. In particular, this provides us with an Eilenberg--Mac~Lane $R$-module spectrum $S$ with $\pi_0 S\cong S_0$ and a bijection of sets
\[\hom_{\mod_R}(R, S)\cong \hom_{\mod_{\pi_{0} R}}(\pi_0 R, S_0).\]
\end{example}

In general though, the answer is more complicated. For a nonexample, consider the periodic real $K$-theory spectrum $\mathrm{KO}$ and the $\pi_\ast \mathrm{KO}$-algebra $S_\ast=\pi_\ast \mathrm{KO} \otimes \mathbb{F}_2$. Using Toda brackets one can show there is no $\mathrm{KO}$-module spectrum $S$ with an isomorphism of $\pi_\ast \mathrm{KO}$-modules $\pi_\ast S\cong S_\ast$; see \cite[Lemma 8.4]{sagave}.\\

This question is also interesting when we consider multiplicative structures. If the spectrum $R$ of Example~\ref{emexample} is an $\E_\infty$-ring spectrum and $\pi_0R\rightarrow S_0$ a map of commutative rings, then $S$ obtains an essentially unique $\E_\infty$-structure such that $R\rightarrow S$ is a map of $\E_\infty$-ring spectra; see \cite[Proposition 7.1.3.18]{ha}. As expected, there are also nonexamples in the multiplicative setting too. Consider the map of rings $\Z\rightarrow \Z(i)$, where the codomain is the ring of Gaussian integers. It is shown in \cite[Proposition 2]{svw} that one cannot construct an $\E_\infty$-ring spectrum $\sph(i)$ lifting (in the sense of \cite[Definition 1]{svw}) the map $\Z\rightarrow \Z(i)$, however, an additive construction is simple -- just take $\sph\vee\sph$. Notice this map $\Z\rightarrow\Z(i)$ is ramified at the prime 2, hence it is not \'{e}tale; see Examples~\ref{absolutelegaloisfailiure} and~\ref{successofinvertingtwo} for more discussion.\\

One solution to the multiplicative problem can be obtained by paraphrasing the work of Baker--Richter using obstruction theory for $\E_\infty$-ring spectra.

\begin{theorem}[{Baker--Richter \cite{br}}]\label{exampleofbr}
Let $A$ be an $\E_\infty$-ring spectrum and $\eta_\ast\colon\pi_\ast A\rightarrow B_\ast$ a map of graded commutative rings recognising $B_\ast$ as a projective $\pi_\ast A$-module. Then there is a homotopy commutative $A$-algebra spectrum $B$ and a map of homotopy commutative ring spectra $\eta\colon A\rightarrow B$, such that $\pi_\ast^e \eta=\eta_\ast$. If in addition, $\eta_\ast$ is \'{e}tale, then $B$ has an $\E_\infty$-structure (unique up to contractible choice) and $\eta$ is a map of $\E_\infty$-ring spectra.\footnote{For the existence of the homotopy commutative $A$-algebra spectrum $B$, one can use the same arguments as in the proof of \cite[Theorem 2.1.1]{br}, as all that is important there is the fact $B_\ast$ is projective over $\pi_\ast A$; see Sections~\ref{additivesection} and~\ref{sectionhomotopycommutative} for more details. For the $\E_\infty$-structure, one can use the same arguments as in the proof of \cite[Proposition 2.2.3]{br}, as the vitally important extra assumption is that $\pi_\ast A\to B_\ast$ is \'{e}tale; see Section~\ref{operads} for more details.}
\end{theorem}

The goal of this article is to explore this question of realisablity and extend Theorem~\ref{exampleofbr} to the setting of \emph{global homotopy theory}, and in this way obtain new global homotopy types. In global homotopy theory, one has global spectra $X$, objects of a stable model category, who naturally have \emph{global homotopy groups}, denoted as $\pi_\ast^G X$ for each compact Lie group $G$. In particular, each global spectrum $X$ has {nonequivariant homotopy groups}, which are simply $\pi_\ast^G X$ when $G=e$ is the trivial group. This concept of a universal equivariant stable homotopy theory has been explored by Bohmann, Greenlees, Lewis, May, and Steinberger; see \cite{b}, \cite{gm}, and \cite{lms}. We will be using the category of orthogonal spectra with the global model structure as defined by Schwede in \cite[Theorem 4.3.18]{s}. This (model) category of global spectra $\spec^\gl$ is symmetric monoidal, so we can speak of monoids (which we call global ring spectra) and commutative monoids (which we call ultra-commutative ring spectra), as well as modules and algebras over these various types of monoids. There also exist intermediary multiplicative structures of global spectra, such as homotopy associative and commutative, $\E_\infty$-global, and $\G_\infty$-ring spectra; see Definitions~\ref{grs},~\ref{einf}, and~\ref{ginfintyguys}, as well as Diagram (\ref{relatingdiagram}) which explains how these concepts relate.\\

To generalise Theorem~\ref{exampleofbr}, one needs to keep in mind that we are not just looking for \textbf{any} realisation of nonequivariant algebraic information by global spectra, but rather realisations over which we understand their global homotopy type. For example, the global spectra $H\mathbf{A}$, $H\mathbf{RU}$, and $H\underline{\Z}$, the \emph{global Eilenberg--Mac~Lane spectra} of the global Burnside ring, global complex representation ring, and constant global functor of $\Z$ (see Remark~\ref{nonexampleforgf} for more details) all have the nonequivariant homotopy type of the Eilenberg--Mac~Lane spectrum $H\Z$, but wildly different global homotopy groups. To overcome this problem, we investigate a condition called \emph{globally flat}; see Definition~\ref{globallyflat}.

\begin{mydef}
Let $R$ be a global ring spectrum and $M$ a left $R$-module spectrum. We say that $M$ is \emph{globally flat} as an $R$-module if a certain natural map
\[\pi_\ast^G R\otimes_{\pi_\ast^e R}\pi_\ast^e M\to \pi_\ast^G M\]
is an isomorphism, for all compact Lie groups $G$. An $R$-algebra is called globally flat if it is globally flat as an $R$-module.
\end{mydef}

Our main theorem then shows that given an ultra-commutative ring spectrum $R$,  certain maps $\pi_\ast^e R\rightarrow S_\ast$ of commutative $\pi_\ast^e R$-algebras can be realised by maps of globally flat $R$-modules $R\rightarrow S$, and that a variety of multiplicative structures can be placed on such an $S$. The following theorem summarises Theorem~\ref{theoremone}, Corollary~\ref{onlyusedforglaoisextensions}, and Theorem~\ref{theoremfive}. 

\begin{theorema}
Let $R$ be an ultra-commutative ring spectrum and $\eta_\ast\colon \pi_\ast^e R\rightarrow S_\ast$ a map of graded commutative rings recognising $S_\ast$ as a projective $\pi_\ast^e R$-module. Then there exists a globally flat homotopy commutative global $R$-algebra $S$, unique to global homotopy. If in addition $\eta_\ast$ is \'{e}tale, then $S$ can be given an $\E_\infty$-global $R$-algebra structure, unique up to contractible choice, lifting the homotopy commutative multiplication. Analogously, if $S_0$ is a $\Q$-algebra, then $S$ has a $\G_\infty$-structure lifting the homotopy commutative multiplication.
\end{theorema}

To prove the above theorem we need to further develop the tools in global homotopy theory a little beyond \cite{s}. In particular, we will relativise some statements made in \cite{s} from the stable global homotopy category $\ho^\gl(\spec)$ to the stable global homotopy category of $R$-modules $\ho^\gl(\mod_R)$, and constantly work with the adjective \emph{globally flat}. As a result, we mimic an array of constructions from classical stable homotopy theory in the setting of global homotopy theory whilst maintaining sufficient control of global homotopy types. For example, one can perform simple localisation constructions, realise Galois extensions of graded rings, and lift nonequivariant spectra from chromatic homotopy theory, all to the global setting. More explicitly, the following is shown as a series of examples in Section~\ref{sectionwithexamples}.\\

\begin{theoremb}
Let $R$ be a fixed ultra-commutative or cofibrant $\E_\infty$-global ring spectrum (see Definition~\ref{einf}), and write $\mathbf{KU}$ and $\MU$ for the \emph{periodic global complex $K$-theory} and \emph{global complex cobordism} spectra; see \cite[Section 6.4]{s} and \cite[Section 6.1]{s}, respectively.
\begin{enumerate}
\item For any (countable) subset $S\subseteq \pi_\ast^e R$, there exists globally flat $\E_\infty$-global $R$-algebra $R[S^{-1}]$ with
\[\pi_\ast^e (R[S^{-1}])\cong (\pi_\ast^e R)[S^{-1}];\]
see Example~\ref{countablesetofelements}. Moreover, for every $R$-module $M$, there exists a globally flat $R$-module $M[S^{-1}]$ with $\pi_\ast^e (M[S^{-1}])\cong (\pi_\ast^e M)[S^{-1}]$; see Example~\ref{localisingmodules}.
\item If a prime $p$ is invertible in $\pi_0^e R$ and the $p^n$th cyclotomic polynomial is irreducible over $\pi_0^e R$, then there exists a globally flat $\E_\infty$-global ring spectrum $R(\zeta)$ realising the map of rings $\pi_0^e R\to (\pi_0^e R)(\zeta)$, where $\zeta$ is a $p^n$th root of unity; see Example~\ref{copingsvw}.
\item If $\pi_\ast^e R\to S_\ast$ is a $G$-Galois extension of (graded) rings for a finite group $G$, then $S_\ast$ is realisable as a globally flat $\E_\infty$-global $R$-algebra $S$, and the $G$-action on $S_\ast$ is realisable by a $G$-action of $\E_\infty$-global $R$-algebras on $S$; see Example~\ref{copingbr}.
\item Every $\pi_\ast^e \mathbf{KU}$-module is realisable by a globally flat $\mathbf{KU}$-module; see Example~\ref{alleasymodulesarerealisable}.
\item For every prime $p$ and every integer $n\geq 0$, there exists a globally flat homotopy associative $\MU_{(p)}$-algebra $\KK(n)$, a \emph{global height $n$ Morava $K$-theory spectrum}, which is nonequivariantly equivalent to the Morava $K$-theory spectrum $\mathrm{K}(n)$; see Example~\ref{moravaktheory}.
\end{enumerate}
The uniqueness of the examples above is also discussed.
\end{theoremb}

Let us now explain the ingredients of this article.

\subsection*{Outline}
In \textsection\ref{backgrounsection}, we recall some of the basic concepts and constructions of global homotopy theory (the details of which can be found inside Schwede's book \cite{s}), in \textsection\ref{homtopabsdo}, we relativise some of this content with respect to a global ring spectrum $R$, and in \textsection\ref{globalyflatrmodules}, globally flat $R$-modules are defined and discussed. The next four sections realise nonequivariant algebraic data in terms of global homotopy theory, first additively in \textsection\ref{additivesection}, multiplicatively up to a single homotopy in \textsection\ref{sectionhomotopycommutative}, multiplicatively up to higher homotopies in \textsection\ref{operads}, and multiplicatively with power operations in \textsection\ref{ginfinity}. In \textsection\ref{additivesection}, we study classical constructions and results (of \cite[Chapter IV]{ekmm}, \cite{wolbert}, and folklore) in the global setting by carefully tracking global flatness, and in \textsection\ref{sectionhomotopycommutative}, we use the ideas of Baker--Richter \cite[Section 2]{br} are also applied to the global homotopy category. In \textsection\ref{operads}, we state and prove some known results about endomorphism operads to help us use the nonequivariant $\E_\infty$-obstruction theory of Goerss--Hopkins \cite{gh04} and Robinson \cite{r}; this section is by far the most technical in this article. In \textsection\ref{ginfinity}, we place $\mathbb{G}_n$-structures (equivalent to certain equivariant norm or multiplicative transfer structures) on certain homotopy commutative global $R$-algebras $S$ when working rationally, essentially as a corollary of \textsection\ref{additivesection} and \textsection\ref{sectionhomotopycommutative}. In \textsection\ref{sectionwithexamples}, we see examples of many of the statements made throughout the rest of the article, and construct new global homotopy types by enriching known nonequivariant and algebraic constructions with controllable global data.

\subsection*{Conventions}
\subsubsection*{Algebraic}
All $G$-representations are finite-dimensional, real, and orthogonal. Homomorphisms of graded rings and graded modules are all degree preserving and graded commutative rings satisfy graded commutativity, $xy=(-1)^{|x||y|}yx$. Given two integers $n,m$ and a graded module $M$, the $m$th level of the shifted module $M[n]$ is $M_{m-n}$, so that for all spectra $X$ we have
\[\pi_m \Sigma^n X\cong \pi_{m-n} X\cong \left((\pi_\ast X)[n]\right)_m.\]

\subsubsection*{Categorical}
All categories are locally small, i.e., the mathematical object $\hom_\C(A,B)$ is always a \emph{set} for each pair of objects $A$ and $B$ in a category $\C$. Let $(\C,\otimes, \1)$ be a symmetric monoidal category. The categories of left and right $R$-modules will be denoted as $\lmod_R$ and $\rmod_R$ respectively. When $R$ is a commutative monoid, $\mod_R$ will denote the category of $R$-modules. An $n$-fold monoidal product over $R$, $M\otimes_R\cdots\otimes_R M$, inside $\mod_R$ will be written as $M^{\otimes n}$. All statements made in this article work equally well for left or right $R$-modules, with the appropriate changes made. 

\subsubsection*{Global homotopical}
The entirety of this article takes place with respect to an arbitrary multiplicative global family $\F$; see \cite[Definition 1.4.1, Proposition 1.4.12(iii)]{s}. This means phrases such as ``global equivalence'', ``global model structure'', ``globally flat'', etc., are all made relative to this global family $\F$. This added flexibility gives us maximum generality, and includes the four most important global families, those of all, finite, and abelian compact Lie groups, as well as the trivial family. This last global family reduces this whole article to statements about the nonequivariant orthogonal spectra of \cite{mmss}. Let $X$ be a global spectrum. For positive integers $n$ we define $\Sigma^n X$ as $X\smp F_{e,0}S^n$ and $\Sigma^{-n}X$ as $X\smp F_{e, \mathbb{R}^n}S^0$ using the notation of \cite[Construction 4.1.23]{s}. The global spectra $F_{e,0}S^n$ and $F_{e,\mathbb{R}^n}S^0$ are globally cofibrant in the model structure of Theorem~\ref{modelstructures}.

\subsubsection*{Model categorical}
Given a topological category $\M$ then $\Map_\M(A,B)$ will denote the mapping space between objects $A$ and $B$ of $\M$. If $\M$ has a topological model structure, the above mapping space does not necessarily have the ``correct homotopy type'' unless $A$ is cofibrant and $B$ is fibrant in $\M$. We will write $[X,Y]_\M = \hom_{\ho(\M)}(X,Y)$. If $R$ is a global ring spectrum (see Definition~\ref{grs}) and $\M=\lmod_R^\gl$ (see Theorem~\ref{globalmodelstructureonrmod}) we will write $[-,-]_\M=[-,-]_R^\gl$. With respect to this model structure the functor $\smp_R\colon \rmod_R^\gl\times \lmod_R^\gl\rightarrow \spec$ is \emph{not} homotopical in either variable and will often be left derived, which can be modelled by taking a cofibrant replacement in either the left, the right, or both variables; see \cite[Theorem 4.3.27]{s}. If $R$ is an ultra-commutative ring spectrum (see Definition~\ref{ucom}) then the homotopy category $\ho^\gl(\mod_R)$ has a symmetric monoidal structure with product $\smpl_R$ and unit $R$, which follows from \cite[Corollary 4.3.29 (ii)]{s}. Homotopy limits and colimits are defined as in \cite[Chapters XI-XII]{bk}, i.e., by deriving the adjoint to the constant diagram functor. For each model category $\M$ that we use, fix cofibrant and fibrant replacement functors $(-)_c$ and $(-)_f$. We assume the reader is familiar with texts on model categories such as \cite{ds}.

\subsubsection*{Topological}
Denote by $\top$ the category of compactly generated weak Hausdorff spaces and continuous maps (see \cite[Appendix A]{s}), which we will call the category of spaces. Let $\ast$ denote the point in $\top$ and write $\top_\ast=\top_{\ast/}$ for the category of based spaces. Given a $G$-representation $V$, then denote by $S^V$ the one-point compactification of $V$, with the $G$-action inherited from $V$. We assume the reader is familiar with the foundations and basics of modern stable homotopy theory from either \cite{ekmm}, \cite{ha}, or \cite{mmss}.

\subsection*{Acknowledgements}
The core of this article is a product of my Masters' thesis at the University of Bonn. I am very happy to thank my Master's advisor Stefan Schwede for suggesting this topic and answering countless questions, as well as my PhD advisor Lennart Meier for many discussions and his valuable insight, and Yuqing Shi for her proof reading. Finally, I greatly appreciate the recommendations and corrections of an anonymous referee, who directly improved the quality of this article.


\section{Background in global homotopy theory}\label{backgrounsection}
Global homotopy theory is the study of spectra with compatible actions of all compact Lie groups $G$ in some \emph{global family} $\F$, a collection of compact Lie groups closed under isomorphisms, closed subgroups, and quotient groups; see \cite[Definition 1.4.1]{s}. We will work with orthogonal spectra as this category, in a certain sense,``contains enough symmetry'' to model global spectra. All of the material in this section can be found in \cite{s} unless otherwise stated.\\

First, let us define $\mathbf{O}$ as the topological category whose objects are real inner product spaces, and whose morphism spaces are defined as
\[\Map_{\mathbf{O}}(V,W)=\mathrm{Th}\left\{(w,\varphi)\in W\times \mathbf{L}(V,W) \enskip |\enskip w\perp \varphi(V)\right\},\]
where $\mathrm{Th}\{\xi\}$ denotes the Thom space of a vector bundle $\xi$, and $\mathbf{L}(V,W)$ the space of linear isometric embeddings from $V$ to $W$. Composition in $\mathbf{O}$ is described in \cite[Construction 3.1.1]{s}. Notice that if $\dim V=\dim W$, then $\mathbf{O}(V,W)$ is homeomorphic $O(V)_+\cong O(W)_+$, the orthogonal groups of $V$ and $W$ with an added basepoint.

\begin{mydef}[{\cite[Example 4.4]{mmss}}]
An \emph{orthogonal spectrum} is a topologically enriched functor $\mathbf{O}\rightarrow \top_\ast$. A map of orthogonal spectra is a natural transformation. Let us denote the category of orthogonal spectra by $\spec$.
\end{mydef}

For us the word spectrum will mean orthogonal spectrum. This category of spectra has a symmetric monoidal structure with product $\smp$ and unit object the sphere spectrum $\sph$; see \cite[Section 21]{mmss} or \cite[Section 3.5]{s}. Following the notation of \cite{ha} we will also write $\oplus$ for the wedge (coproduct) of spectra.\\

We now make a crucial observation. Let $X$ be a spectrum, $G$ be a compact Lie group in $\F$, and $V$ be any $G$-representation. By considering $V$ as a real inner product space, we obtain a based space $X(V)$, which by functorality of $X$ has a $G$-action: 
\[G\longrightarrow O(V)_+\cong\mathbf{O}(V,V)\overset{X}{\longrightarrow} \Map_{\top_\ast}(X(V), X(V)).\]
This is how the category of orthogonal spectra encodes the representation theory of all compact Lie groups $G$, and in a certain sense ``contains enough symmetry''.

\begin{mydef}[{\cite[Definition III.3.2]{mm}}]\label{equivarianthomotopygroups}
Let $X$ be a spectrum and $G$ a compact Lie group. We define the \emph{zeroth $G$-homotopy group of $X$} as the colimit
\[\pi_0^G X=\underset{V\in s(\U_G)}{\colim} [S^V,X(V)]^G_\ast,\]
where $[-,-]^G_\ast$ denotes homotopy classes of continuous equivariant maps of based $G$-spaces, $s(\U_G)$ denotes the poset of finite $G$-subrepresentations of the complete $G$-universe; see \cite[Definition 1.1.12]{s}, and the maps in the colimit are defined by the composition
\[[S^V,X(V)]^G_\ast\xrightarrow{S^U\smp -} [S^{U\oplus V},S^U\smp X(V)]^G_\ast\xrightarrow{(\sigma_{U,V}\colon S^U\smp X(V)\to X(U\oplus V))_\ast} [S^{U\oplus V},X(U\oplus V)]^G_\ast,\]
where the latter map is postcomposition with a certain structure morphism of $X$; see \cite[p.232]{s}.
\end{mydef}

The above definition does not depend on the chosen complete $G$-universe by a cofinality argument; see \cite[Remark V.1.10]{mm}. To define homotopy groups $\pi_q^G X$ for $q\neq 0$ we either smash the domain $S^V$ with $S^q$ on the right for $q>0$, or shift the codomain $V$ on the right by $\R^{-q}$ for $q<0$; see \cite[(3.1.11)]{s}. These sets $\pi_q^G X$ have a natural abelian group structure for all compact Lie groups $G$ and all integers $q$; see \cite[p.233]{s}. Write $\pi_\ast^G X$ for the graded abelian group $\bigoplus_{q\in \Z}\pi_q^G X$. There is a wealth of structure between $\pi_q^G X$ and $\pi_q^K X$ for two compact Lie groups $G$ and $K$ in $\F$. For every continuous homomorphism of compact Lie groups $\al\colon K\rightarrow G$ there is a restriction map $\al^\ast\colon \pi_q^G X\rightarrow \pi_q^K X$, which is constructed by pulling $G$-actions back to $K$-actions; see \cite[Construction 3.1.15]{s}. For each closed inclusion of compact Lie groups $H\leq G$ there is a transfer map $\tr_H^G\colon \pi_q^H X\rightarrow \pi_q^G X$, which is defined using a Thom-Pontryagin construction; see \cite[Section 3.2]{s}. These two families of maps generate the set of natural transformations from $\pi_0^G$ to $\pi_0^K$ as functors from $\spec$ to $\ab$, which are the natural operations on global homotopy groups; see \cite[Proposition 4.2.5 \& Theorem 4.2.6]{s}.

\begin{mydef}[{\cite[Construction 4.2.1, Definition 4.2.2]{s}}]
Let $\mathbf{A}$ be the pre-additive \emph{global Burnside category}, whose objects are compact Lie groups inside $\F$ and morphism groups are defined by
\[\hom_{\mathbf{A}}(G,K)=\nat(\pi_0^G, \pi_0^K).\]
A \emph{global functor} is an additive functor from $\mathbf{A}$ to the category of abelian groups. Let $\gf$ denote the category of global functors and natural transformations.
\end{mydef}

By definition the assignment $X\mapsto \{\pi_q^G X\}_{G\in\F}=\epi_q X$ constitutes a global functor for any spectrum $X$ and any integer $q$. We define a \emph{graded global functor} to be a collection of global functors $\{F_q\}_{q\in\mathbb{Z}}$. For any global spectrum $X$ we write $\epi_\ast X$ for the graded global functor $\bigoplus_{q\in\mathbb{Z}}\epi_q X$.\\

Suppose we have a map of orthogonal spectra $f\colon X\rightarrow Y$, then by Definition~\ref{equivarianthomotopygroups} we see the construction of equivariant homotopy groups is functorial. We obtain an array of induced maps, for all compact Lie groups $G$ in $\F$ and integers $q$, which we all call $f_\ast$,
\[f_\ast\colon \pi_q^G X\longrightarrow \pi_q^G Y, \qquad f_\ast\colon \pi_\ast^G X\longrightarrow \pi_\ast^G Y,\qquad f_\ast\colon \epi_q X\longrightarrow \epi_q Y,\qquad f_\ast\colon \epi_\ast X\longrightarrow \epi_\ast Y.\]

\begin{mydef}[{\cite[Definition 4.1.3]{s}}]
Let $f\colon X\rightarrow Y$ be a map of orthogonal spectra. We say $f$ is a \emph{global equivalence} if the induced map $f_\ast\colon \epi_\ast X\rightarrow \epi_\ast Y$ is an isomorphism.
\end{mydef}

A theorem of Schwede says the global equivalences are part of a model structure on $\spec$.

\begin{theorem}[{\cite[Theorem 4.3.17]{s}}]\label{modelstructures}
There exists a topological stable model structure on $\spec$, the \emph{global model structure}, whose weak equivalences are global equivalences and fibrant objects the global-$\Omega$-spectra; see \cite[Definition 4.3.14]{s}.
\end{theorem}

Denote by $\spec^\gl$ the category of orthogonal spectra with the global model structure of Theorem~\ref{modelstructures}. We remind the reader that for us the phrases ``global equivalence'' and ``global-$\Omega$-spectra'' are relative to an ambient global family $\F$. We will write $\ho^\gl(\spec)$ for the homotopy category of $\spec^\gl$. When the ambient global family is trivial (when $\F=\{e\}$ contains only the trivial group) then $\spec^\gl$ will be written as $\spec^e$, which is \emph{equal} to the stable model category of orthogonal spectra defined in \cite[Theorem 9.2]{mmss}.

\begin{remark}\label{modelstructurestwo}
In particular, by \cite[Definition 4.3.14]{s} we see a global equivalence of spectra is a nonequivariant equivalence, and a global fibration is a nonequivariant fibration. This also implies a nonequivariant cofibration is a global cofibration by standard model categorical lifting properties; see \cite[Proposition 3.13]{ds}.
\end{remark}

In this article we would like to study global homotopy theory relative to a ring spectrum $R$. There are many different types of ring spectra one can talk about, with various levels of multiplicative structure. Let us first make the purely categorical definitions.

\begin{mydef}[{\cite[Definition 3.5.15]{s}}]\label{grs}
A \emph{global ring spectrum} is a monoid object of $\spec^\gl$. A \emph{homotopy associative (resp. commutative) global ring spectrum} is an associative (resp. commutative) monoid object of $\ho^\gl(\spec)$.
\end{mydef}

Let us recall some operadic definitions.

\begin{mydef}
A \emph{topological monoidal model category} is a topological model category (see \cite[Definition 5.12]{mmss}) endowed with a closed symmetric monoidal structure which satisfies the \emph{pushout product axiom} of \cite[Definition 3.1]{ss}.
\end{mydef}

Suppose $(\M, \otimes, \1)$ is a topological monoidal model category, then for any object $X$ of $\M$ the $n$th level of the \emph{endomorphism operad of $X$} is defined as the mapping space
\[(\opend_\M X)_n=\Map_\M(X^{\otimes n}, X)\]
with the tautological $\Sigma_n$-action from $X^{\otimes n}$. Let $\O$ be a topological operad. An \emph{$\O$-algebra in $\M$} is a map of topological operads $\ga\colon \O\rightarrow \opend_\M(X)$, which is only homotopically well-defined if $X$ is a bifibrant object of $\M$. The category of topological operads has a model structure, with weak equivalences (resp. fibrations) gives by levelwise topological weak equivalences (resp. fibrations); see \cite[Example 3.3.2]{bm}. An \emph{$\E_\infty$-operad} is a $\Sigma$-cofibrant replacement of the terminal (commutative) operad, and an \emph{$\E_\infty$-object in $\M$} is an $\O$-algebra in $\M$ for any $\E_\infty$-operad $\O$; see \cite[Section 1]{bm}. For a proof that the definition of an $\E_\infty$-object is independent of the chosen $\E_\infty$-operad (see \cite[Section 4]{bm}), but for consistency, let us fix a topological $\E_\infty$-operad $\O$.

\begin{mydef}\label{einf}
An \emph{$\E_\infty$-global ring spectrum} is an $\E_\infty$-object of $\spec^\gl$. 
\end{mydef}

By \cite[Theorem 4.4]{spit} (or \cite[Example 4.6.4]{bm}), the category of $\E_\infty$-global ring spectra, denoted as $\calg^\gl$, has an induced model structure from $\spec^\gl$ (as the latter satisfies the monoid axiom by \cite[Proposition 4.3.28]{s}), so weak equivalences (resp. fibrations) are given by global weak equivalences (resp. global fibrations) in $\spec^\gl$.\\

The same holds for the trivial global family $\F=\{e\}$, and we denote by $\calg^e$ the model category of nonequivariant $\E_\infty$-ring spectra, called $\E_\infty$-rings. Moreover, with these definitions, we see the identity $\calg^\gl\to \calg^e$ is a right Quillen functor (with left adjoint also given by the identity); this is further justified by \cite[Theorem 2.14]{togetmodulesovereinfty}.\\

Let us warn the reader that an $\E_\infty$-global ring spectrum is not in general globally equivalent to a strictly commutative orthogonal spectrum (unless the global family $\F$ is trivial). There is a tactable difference between these two notions of commutativity in equivariant and global homotopy theory: multiplicative norms and power operations; see \cite{blumberghill} and \cite[Section 5]{s}, respectively.

\begin{mydef}[{\cite[Definition 5.1.1]{s}}]\label{ucom}
An \emph{ultra-commutative ring spectrum} is a commutative monoid of $\spec^\gl$.
\end{mydef}

The sphere spectrum $\sph$, the Thom spectra $\mathbf{MO}$ and $\MU$, and the connective global $K$-theory spectrum $\mathbf{ku}$ are all ultra-commutative ring spectra; see \cite[p.303]{s}, \cite[Section 6.1]{s}, and \cite[Construction 6.3.9]{s}, respectively. The $\E_\infty$-global ring spectrum $\mathbf{mO}$ of \cite[p.303]{s} is \emph{not} ultra-commutative, as demonstrated by a lack of power operations. Let $p$ be a prime greater than 3, then the Moore spectra $\sph/p$ (the cofibres of multiplication by $p\colon\sph\rightarrow\sph$) are examples of homotopy commutative but not $\E_\infty$-global or even simply global ring spectra; see \cite[Example 3.3]{vig}. There is also a concept of a homotopy commutative global spectrum with power operations, called \emph{$\G_\infty$-ring spectra}, which mimic the nonequivariant $\mathbb{H}_\infty$-ring spectra of \cite{brunermaymcluresteinberger}. These are not $\E_\infty$-global ring spectra by lifting the nonequivariant example of \cite{einftyandhingty} into global homotopy theory; see \cite[Example 3.46]{m}.\\

In summary, we have the following diagram of implications between adjectives of global spectra.\footnote{An $\E_\infty$-global ring spectrum is an $\mathbb{A}_\infty$-global ring spectrum (using the definition of an $\mathbb{A}_\infty$-object from \cite[Remark 4.6]{bm}) as a cofibrant replacement of the unique map from the associative operad to the commutative operad implies all $\E_\infty$-algebras are $\mathbb{A}_\infty$-algebras. An application of \cite[Remark 4.6]{bm} in $\spec^\gl$ shows the model categories of $\mathbb{A}_\infty$-global ring spectra and global ring spectra are Quillen equivalent. In particular, there is also an arrow in Diagram (\ref{relatingdiagram}) from $\E_\infty$-global ring spectra to global ring spectra, but we will not use this fact.}

\begin{equation}\label{relatingdiagram}\begin{tikzcd}
{\mbox{ultra-commutative ring spectrum}}\ar[r, Rightarrow]\ar[dd, Rightarrow]\ar[rd, Rightarrow]	&	{\mbox{global ring spectrum}}	\\
	&	{\G_\infty\mbox{-ring spectrum}}\ar[d, Rightarrow]	\\
{\E_\infty\mbox{-global ring spectrum}}\ar[r, Rightarrow]	&	{\mbox{homotopy commutative global ring spectrum}}	
\end{tikzcd}.\end{equation}

For a global ring spectrum $R$ we have a categories of left and right $R$-modules, which obtain global model structures through the extension of scalars adjunction

\begin{equation}\label{adjectunoi}\Map_{\lmod_R^\gl}(R\smp X,M)\cong \Map_{\spec^\gl}(X, M).\end{equation}

\begin{theorem}[{\cite[Corollary 4.3.29]{s}}]\label{globalmodelstructureonrmod}
Let $R$ be a global ring spectrum. There are topological model structures on $\lmod_R$ and $\rmod_R$ whose weak equivalences (resp.\ fibrations) are the weak equivalences (resp.\ fibrations) of $\spec^\gl$. Moreover, if $R$ is an ultra-commutative ring spectrum, then $\mod_R$ is a monoidal model category with respect to $\smp_R$.
\end{theorem}

Denote by $\lmod_R^\gl$ and $\rmod_R^\gl$ the topological monoidal model categories given above. In particular, when our ambiant global family $\F$ is the trivial global family, we will write $\lmod_R^e$ and $\rmod_R^e$, which are \emph{equal} to the nonequivariant model categories of left and right $R$-module orthogonal spectra of \cite[Theorem 12.1]{mmss}. Taking $R=\sph$, we see $\spec^\gl$ is also a topological monoidal model category.\\

Definitions~\ref{grs}, ~\ref{einf}, and~\ref{ucom} can all be relativised (by taking categories under $R$) to define global $R$-algebras, homotopy commutative $R$-algebras, $\E_\infty$-global $R$-algebras, and ultra-commutative $R$-algebras, respectively. In particular, the category of $\E_\infty$-global $R$-algebras will be given a model structure by considering it as the category of $\E_\infty$-global ring spectra under a fixed $R$; see Definition~\ref{einf} and \cite[Remark 3.10]{ds}. With this definition, the identity $\calg^\gl_R\to \calg_R^e$ is a right Quillen functor, with left adjoint the identity too.


\section{Homotopy theory over a global ring spectrum}\label{homtopabsdo}
In \cite{s}, the foundations of global homotopy theory were mostly established over the global sphere spectrum. In this section we will extend some results of \cite[Section 4]{s} to statements over an arbitrary global ring spectrum $R$.

\begin{prop}[{\cite[Proposition 4.3.22(i), Theorem 4.4.3]{s}}]\label{compactlygenerated}
Let $R$ be a global ring spectrum. The triangulated category $\ho^\gl(\lmod_R)$ is compactly generated and has coproducts indexed on arbitrary sets.
\end{prop}

\begin{proof}
Using \cite[Theorem 4.4.3]{s} and the fact (\ref{adjectunoi}) is a Quillen adjunction shows the set of $R$-modules $\{R\smp \Sigma_+^\infty B_\gl G\}_{G\in\F}$ is a set of compact weak generators of $\ho^\gl(\lmod_R)$; see \cite[Definition 1.1.27, Construction 4.1.7]{s}. The statement about coproducts follows by the same argument from \cite[Proposition 4.3.22(i)]{s}, as coproducts in $\ho^\gl(\lmod_R)$ can be modelled by a wedge of bifibrant objects in $\lmod_R$.
\end{proof}

\begin{construction}\label{constr}
The proof of \cite[Theorem 4.4.3]{s} uses the fact the spectra $\Sigma_+^{\infty+q} B_\gl G$ represent the functors $\pi_q^G$ from $\ho^\gl(\spec)\rightarrow \ab$. If $R$ is a global ring spectrum, then the fact Adjunction~\ref{adjectunoi} is a Quillen adjunction with respect to the model structures of Theorems~\ref{modelstructures} and~\ref{globalmodelstructureonrmod} means the left $R$-module $R\otimes \Sigma_+^{\infty+q} B_\gl G$ represents the functor
\[\pi_q^G\colon\ho^\gl(\lmod_R)\rightarrow \ab.\]
This means that given a fibrant left $R$-module $M$ and an element $x\in \pi_q^e M$, then we can represent $x$ by a map of left $R$-modules $R\otimes \Sigma_+^{\infty+q}B_\gl e=\Sigma^q R\rightarrow M$.
\end{construction}

\begin{prop}[{\cite[Theorem 4.5.1]{s}}]\label{needleftadjunction}
Let $R$ be a global ring spectrum. Then the identity functor $\id\colon \lmod_R^\gl\rightarrow \lmod_R^e$ is a right Quillen functor, whose derived left adjoint $L\colon \ho(\lmod_R)\rightarrow \ho^\gl(\lmod_R)$ is fully faithful.
\end{prop}

\begin{proof}
From the definitions of the model structures on $\lmod_R^\gl$ and $\lmod_R^e$, we see the identity $\id\colon \lmod_R^\gl\rightarrow \lmod_R^e$ is a right Quillen functor with $\id\colon \lmod_R^e\rightarrow \lmod_R^\gl$ the associated left Quillen functor. The right Quillen functor takes all global equivalences to weak equivalences so it need not be derived to induce a functor $U$ on homotopy categories. The unit $\eta_M\colon M\rightarrow ULM$ of the derived adjunction
\[\left[LM, N\right]_R^\gl\cong \left[M, UN\right]_R^e,\]
is then an isomorphism for all objects $M$ of $\ho(\lmod_R)$. Hence $L$ is fully faithful.
\end{proof}

\begin{mydef}[{\cite[Definition 4.5.6]{s}}]\label{leftindcued}
Let $R$ be a global ring spectrum. We say a left $R$-module $M$ is \emph{left induced} if $M$ is in the essential image of the functor $L$.
\end{mydef}

\begin{remark}[{\cite[Remark 4.5.3]{s}}]\label{cofibrantobjectsgivege}
One can calculate the value of $L$ on an $R$-module $M$ by taking a nonequivariant cofibrant replacement $M_c\rightarrow M$ of $M$. The global homotopy type of $M_c$ is then well-defined. Indeed, as $\id\colon \lmod_R^e\rightarrow\lmod_R^\gl$ is a left Quillen functor, then nonequivariant acyclic cofibrations are sent to global acyclic cofibrations, and by Ken Brown's lemma (see \cite[Lemma 9.9]{ds}) we see nonequivariant weak equivalences between nonequivariant cofibrant objects are in fact global equivalences. In particular, we see that the derived adjunction counit $\epsilon_M\colon LUM\to M$ can be modelled by taking a nonequivariant cofibrant replacement of $M$.
\end{remark}

This remark implies the following alternative characterisation of left induced modules.

\begin{cor}\label{somethingsimplethanksref}
Let $R$ be a global ring spectrum and $M$ a left $R$-module. Then the following are equivalent.
\begin{enumerate}
\item The left $R$-module $M$ is left induced.
\item The derived adjunction counit $\epsilon_M\colon LUM\to M$ is an isomorphism in $\ho^\gl(\lmod_R)$.
\item A (and hence every) nonequivariant cofibrant replacement $M_c\to M$ of $M$ in $\lmod_R^e$ is in fact a global equivalence.
\end{enumerate}
\end{cor}

The same statement holds for right $R$-modules, \emph{mutatis mutandis}. For use in this proof, let $e$-$\mathcal{P}$ (resp.\ $\gl$-$\mathcal{P}$) refer to a model categorical property $\mathcal{P}$ inside $\mod_R^e$ (resp.\ inside $\mod_R^\gl$). We will also use Remark~\ref{modelstructurestwo} without mention.

\begin{proof}
Without loss of generality $M$ is $\gl$-bifibrant. By Remark~\ref{cofibrantobjectsgivege}, parts 2 and 3 are equivalent, and part 2 implies part 1 by definition. To see part 1 implies part 3, suppose that $M$ is left induced, then by Remark~\ref{cofibrantobjectsgivege} there exists an $e$-cofibrant $R$-module $N_c$ and an isomorphism $N_c\to M$ in $\ho^\gl(\lmod_R)$. From our (co)fibrancy assumptions, this lifts to a strict map $f\colon N_c\to M$ in $\lmod_R^e$, and factors through an $e$-cofibrant replacement of $M$,
\[N_c\xrightarrow{g} M_c\xrightarrow{h} M,\]
as $N_c$ is $e$-cofibrant and $h$ is an $e$-acyclic fibration. The map $f$ is a $\gl$-equivalence by assumption, by Remark~\ref{cofibrantobjectsgivege} the $e$-equivalence $g$ is also a $\gl$-equivalence, hence $h$ is a $\gl$-equivalence.
\end{proof}


\section{Globally flat $R$-modules}\label{globalyflatrmodules}
Studying the left induced left $R$-modules of Definition~\ref{leftindcued} is one way to safely pass from nonequivariant to global information. However, it is not as tactable as one might like, which leads us to the following.

\begin{mydef}\label{globallyflat}
Let $R$ be a global ring spectrum and $M$ a left $R$-module. We say $M$ is \emph{globally flat} if for all $G$ inside $\F$ the canonical $\pi_\ast^G R$-module morphisms
\[\Lambda_M^G\colon \pi_\ast^G R\underset{\pi_\ast^e R}{\otimes} \pi_\ast^e M\longrightarrow \pi_\ast^G M,\qquad r\otimes m\longmapsto r\cdot p^\ast_G(m)\]
are isomorphisms, where $p_G\colon G\rightarrow e$ is the unique map. An $R$-algebra (of any kind) is globally flat if the underlying $R$-module is.
\end{mydef}

We will see some examples of $R$-modules in Proposition~\ref{naturalnessandstuff} which are globally flat, and there are also natural nonexamples.

\begin{remark}\label{nonexampleforgf}
Consider a \emph{global Eilenberg--Mac~Lane spectrum} (see \cite[Remark 4.4.12]{s}), which is a global spectrum $HF$ associated to a global functor $F$, defined uniquely up to isomorphism in $\ho^\gl(\spec)$ by the requirement that $\epi_0 HF\cong F$ and $\epi_q HF=0$ for all $q\neq 0$. Let us also consider the global functors $\burn$ and $\mathbf{RU}$, which are defined such that for a finite group $G$ the group $\burn(G)$ is the \emph{Burnside ring of finite $G$-sets} and $\mathbf{RU}(G)$ is the \emph{complex representation ring of $G$}; see \cite[Example 4.2.8]{s}. We claim that $H\mathbf{RU}$ could never be globally flat over $H\burn$ (so long as $\F$ is not trivial) as the Burnside ring $A(G)$ and the complex representation ring $RU(G)$ are not isomorphic as abelian groups for all compact Lie groups $G$, the smallest example being $G=C_3$. The same goes for the global Eilenberg--Mac~Lane spectrum of the constant global functor at $\Z$ over $H\mathbf{A}$; see \cite[Example 4.2.8]{s}.
\end{remark}

\begin{remark}\label{naturality}
Notice that for each $G$ in $\F$, the map $\Lambda_M^G$ above is the image of the map $p_G^\ast\colon \pi_\ast^e M\rightarrow \pi_\ast^G M$ under the extension of scalars adjunction induced by $p_G^\ast\colon \pi_\ast^e R\rightarrow \pi_\ast^G R$
\[\hom_{\lmod_{\pi_\ast^G R}}\left(\pi_\ast^G R\underset{\pi_\ast^e R}{\otimes} \pi_\ast^e M, \pi_\ast^G M\right)\cong \hom_{\lmod_{\pi_\ast^e R}}\left(\pi_\ast^e M,\pi_\ast^G M\right).\]
For a global $R$-algebra $S$, as the map $p_G^\ast\colon \pi_\ast^e S\rightarrow \pi_\ast^G S$ is a multiplicative map, then using an extension of scalars adjunction for graded $\pi_\ast^e R$-algebras we see $\Lambda_S^G$ is also multiplicative in this case. Let us summarise some more properties of these maps below.
\end{remark}

\begin{prop}\label{naturalnessandstuff}
Let $R$ be a global ring spectrum. Then for all $G$ in $\F$ the maps $\Lambda_M^G$ are natural in the $R$-module variable $M$, and form a morphism of graded global functors
\[\Lambda_M\colon \epi_\ast R \underset{\pi_\ast^e R}{\otimes} \pi_\ast^e M\longrightarrow \epi_\ast M.\]
Moreover, if $\lmod_R^\Lambda$ denotes the full subcategory of $\lmod_R^\gl$ spanned by the globally flat $R$-modules, then $\lmod_R^\Lambda$ is closed under arbitrary suspensions, wedges, filtered homotopy colimits, and contains $R$.
\end{prop}

\begin{proof}
Defining $\Lambda_M^G$ using the extension of scalars adjunction from Remark~\ref{naturality} shows the naturality in $M$. For naturality in the compact Lie group variable we need to show these maps commute with restrictions and transfers, as \cite[Proposition 4.2.5 \& Theorem 4.2.6]{s} imply these maps form a $\Z$-basis of $\hom_\mathbf{A}(G,K)$, for any $G,K$ in $\F$. Fix some $R$-module $M$, and let $f\colon K\rightarrow G$ be any morphism of compact Lie groups in $\F$. The compatibility of these maps with restrictions then follows from the equalities
\[\left(\Lambda^K_M\circ (f^\ast\otimes \id)\right)\left(r\otimes m\right)=f^\ast r\cdot p^\ast_K m=f^\ast r\cdot f^\ast p^\ast_G m=f^\ast(r\cdot p^\ast_G(m))=\left(f^\ast\circ\Lambda^G_M\right)\left(r\otimes m\right).\]
The second equality comes from the equality $p_K= p_G\circ f$ of group homomorphisms, and the third equality from the that fact restriction maps are $\pi_\ast^G R$-module homomorphisms. For the transfers, let $H$ be a closed subgroup of a compact Lie group $G$ inside $\F$, then we obtain the following equalities,
\[\left(\Lambda_M^G\circ (\tr_H^G\otimes\id)\right)\left(r\otimes m\right)=\tr_H^Gr\cdot p^\ast_G m=\tr_H^G(r\cdot \res^G_H(p_G^\ast m))\]
\[=\tr_H^G(r\cdot p_H^\ast m)=\left(\tr_H^G\circ \Lambda_M^H\right)\left(r\otimes m\right).\]
The second equality is a consequence of Frobenius reciprocity (see \cite[Corollary 3.5.17(v)]{s}) and the third equality from the equality $p_H=p_G\circ i$ as group homomorphisms. This shows the maps $\Lambda_M^G$ are natural in $G$, hence $\Lambda_M$ is a morphism of graded global functors.\\

For the ``moreover statement'', notice $R$ is in $\lmod_R^\Lambda$ as for all $G$ in $\F$ the map $\Lambda_R^G$ is the canonical isomorphism
\[\pi_\ast^G R\underset{\pi_\ast^e R}{\otimes} \pi_\ast^e R\overset{\cong}{\longrightarrow} \pi_\ast^G R.\]
When $M$ is an object of $\lmod_R^\Lambda$, then for all integers $n$, $\Sigma^n M$ is in $\lmod_R^\Lambda$ from the natural isomorphisms
\[\pi_\ast^G R\underset{\pi_\ast^e R}{\otimes}\pi_\ast^e \Sigma^n M\cong \left(\pi_\ast^G R\underset{\pi_\ast^e R}{\otimes}\pi_\ast^e M\right)[n]\cong  (\pi_\ast^G M)[n]\cong \pi_\ast^G \Sigma^n M.\]
If $M_i$ are objects of $\lmod_R^\Lambda$ for all $i$ in some indexing set $I$, then $\bigvee_{i\in I} M_i$ is also in $\lmod_R^\Lambda$ from the natural isomorphisms
\[\pi_\ast^G R\underset{\pi_\ast^e R}{\otimes}\pi_\ast^e \left( \bigvee_{i\in I}M_i\right)\cong \pi_\ast^G R\underset{\pi_\ast^e R}{\otimes}\bigoplus_{i\in I}\pi_\ast^e M_i\cong \bigoplus_{i\in I}\pi_\ast^G M_i\cong \pi_\ast^G \left( \bigvee_{i\in I} M_i \right).\]
Finally, if we have a filtered system of left $R$-modules $M_i$ inside $\lmod_R^\Lambda$, then $\hocolim_i M_i$ is in $\lmod_R^\Lambda$ from the natural isomorphisms
\[\pi_\ast^G R\underset{\pi_\ast^e R}{\otimes}\pi_\ast^e (\underset{i}{\hocolim} \,M_i)\cong \pi_\ast^G R\underset{\pi_\ast^e R}{\otimes}\underset{i}{\colim} (\pi_\ast^e M_i)\cong \underset{i}{\colim}(\pi_\ast^G M_i)\cong \pi_\ast^G (\underset{i}{\hocolim} \,M_i) .\qedhere\]
\end{proof}

\begin{remark}\label{easyconsequence}
One consequence of Definition~\ref{globallyflat} and the naturality of these maps $\Lambda_M^G$ in $M$ is the following simple observation. Let $f\colon M\rightarrow N$ be a map of globally flat $R$-modules. Then $f$ is a global equivalence if and only if $\pi_\ast^e f$ is an equivalence. The ``only if'' direction follows from Definition~\ref{globallyflat}, and the converse is a consequence of the following naturality diagram of $\pi_\ast^G R$-modules,
\[\begin{tikzcd}
{\pi_\ast^G R\underset{\pi_\ast^e R}{\otimes} \pi_\ast^e M}\ar[rr, "{\id\otimes \pi_\ast^e f}"]\ar[d, "\cong", swap]\ar[d, "{\Lambda_M^G}"]\ar[d, "\cong", swap]	&&	{\pi_\ast^G R\underset{\pi_\ast^e R}{\otimes}\pi_\ast^e N}\ar[d, "{\Lambda_{N}^G}"]\ar[d, "\cong", swap]\\
{\pi_\ast^G M}\ar[rr, "{\pi_\ast^G f}"]																								&&	{\pi_\ast^G N}
\end{tikzcd}.\]
Notice how this resembles Remark~\ref{cofibrantobjectsgivege}, in that the global homotopy type of both the classes of left induced and of globally flat $R$-modules are controlled by nonequivariant information.
\end{remark}


\section{Realising algebra with $R$-modules}\label{additivesection}
Our first step in realising algebra in global homotopy theory is additive, i.e., as $R$-modules. 

\begin{prop}\label{extensiontoflatness}
Let $R$ be a global ring spectrum and $M_\ast$ a projective left $\pi_\ast^e R$-module. There is a globally flat left $R$-module $M$ and an isomorphism $\phi_g\colon \pi_\ast^e M\cong M_\ast$ of left $\pi_\ast^e R$-modules.
\end{prop}

\begin{proof}
The projectivity condition means there is an idempotent morphism of $\pi_\ast^e R$-modules
\[f\colon \bigoplus_{i\in I} (\pi_\ast^e R)[n_i] \longrightarrow \bigoplus_{i\in I} (\pi_\ast^e R)[n_i],\]
for some indexing set $I$ and $n_i\in \Z$, and a $\pi_\ast^e R$-module isomorphism $\phi_f\colon \im(f)\rightarrow \pi_\ast^e M$. Define $F$ as a fibrant replacement (in $\lmod_R^\gl$) of $\bigvee_{i\in I} \Sigma^{n_i}R$. By construction $\pi_\ast^e F\cong \bigoplus_{i\in I} (\pi_\ast^e R)[n_i]$. We can construct a map of left $R$-modules $g\colon F\rightarrow F$ such that $\pi_\ast^e g=f$ by Construction~\ref{constr}. This implies $g$ is idempotent in $\ho^\gl(\lmod_R)$. Proposition~\ref{compactlygenerated} allows us to use \cite[Proposition 1.6.8]{amnon} with respect to the idempotent map $g\colon F\rightarrow F$, which gives us a commutative diagram in $\ho^\gl(\lmod_R)$,
\begin{equation}\label{retract}\begin{tikzcd}
{F}\ar[r]\ar[rr, "g", bend right = 25, swap]	&	{\widetilde{M}}\ar[r]\ar[rr, "\id", bend left = 25]	&	{F}\ar[r]	&	{\widetilde{M}}
\end{tikzcd},\end{equation}
where $\widetilde{M}$ is the homotopy colimit of $F\overset{g}{\rightarrow} F\overset{g}{\rightarrow}\cdots$. As $\pi_\ast^e g=f$ we see
\[\pi_\ast^e \widetilde{M}= \pi_\ast^e \hocolim\left(	F\overset{g}{\longrightarrow}F\overset{g}{\longrightarrow}\cdots	\right)\cong \colim\left(\pi_\ast^eF\overset{f}{\longrightarrow}\pi_\ast^eF\overset{f}{\longrightarrow}\cdots\right)\cong \im(f)\cong \pi_\ast^e M,\]
using the fact $f^2=f$. Proposition~\ref{naturalnessandstuff} shows $\widetilde{M}$ is globally flat.
\end{proof}

Next we consider realising morphisms of $\pi_\ast^e R$-modules by morphisms of $R$-modules.

\begin{prop}\label{proptwo}
Let $R$ be a global ring spectrum and $M$ a globally flat left $R$-module such that $\pi_\ast^e M$ is projective as a left $\pi_\ast^e R$-module. Then for all left $R$-modules $N$, the functor $\pi_\ast^e$ induces an isomorphism of abelian groups
\[\left[M,N\right]_R^\gl\overset{\pi_\ast^e}{\longrightarrow} \hom_{\lmod_{\pi_\ast^eR}}(\pi_\ast^e M, \pi_\ast^e N).\]
\end{prop}

\begin{proof}
First let us assume $M$ is a wedge of suspensions of $R$, so
\[M=\bigvee_{i\in I}\Sigma^{n_i}R.\]
Suppose $N$ is an arbitrary left $R$-module and consider
\begin{equation}\label{letussee}\left[R,N\right]_R^\gl\overset{\pi_\ast^e }{\longrightarrow}\hom_{\lmod_{\pi_\ast^e R}}(\pi_\ast^e R, \pi_\ast^e N).\end{equation}
The above map is a bijection as both of the above sets are canonically in bijection with $\pi_0 N$, the left by representability (see Construction~\ref{constr}), and the right by elementary algebra. To extend this observation we consider the following diagram of abelian groups,
\begin{equation}\label{everythinginsniceandfree}\begin{tikzcd}
{\left[\bigvee_{i\in I}\Sigma^{n_i}R,N\right]_R^\gl}\ar[rr, "\pi_\ast^e"]\ar[d, "\cong"]	&&	{\hom_{\lmod_{\pi_\ast^e R}}\left(\bigoplus_{i\in I} (\pi_\ast^e R)[n_i],\pi_\ast^e N\right)}\ar[d, "\cong"]\\
{\prod_{i\in I}\left[\Sigma^{n_i}R,N\right]_R^\gl}\ar[d, "\cong"]\ar[rr, "\prod_{i\in I}\pi_\ast^e"]				&&	{\prod_{i\in I} \hom_{\lmod_{\pi_\ast^e R}}\left((\pi_\ast^e R)[n_i],\pi_\ast^e N\right)}\ar[d, "\cong"]\\
{\prod_{i\in I}\left[R,\Sigma^{-n_i}N\right]_R^\gl}\ar[rr, "\pi_\ast^e"]\ar[rr, "\cong", swap]	&&	{\prod_{i\in I} \hom_{\lmod_{\pi_\ast^e R}}\left(\pi_\ast^e R,(\pi_\ast^e N)[-n_i]\right)}
\end{tikzcd}.\end{equation}
The vertical isomorphisms come from the universal property of coproducts, or properties of shifts. The naturality of these maps give us the commutativity of the above diagram. The lower-horizontal map is a product of the isomorphism (\ref{letussee}) and the quick calculation $\pi_\ast^e \Sigma^{-n_i}N\cong (\pi_\ast^e N)[-n_i]$. This gives us our desired result in the case when $M$ is a wedge of suspensions of $R$.\\

Consider now a globally flat left $R$-module $M$ such that $\pi_\ast^eM$ is projective over $\pi_\ast^e R$. By Proposition~\ref{extensiontoflatness} we have a left $R$-module $\widetilde{M}$ which realises $\pi_\ast^e M$. Using the same notation as in Proposition~\ref{extensiontoflatness}, we see by (\ref{retract}) $\widetilde{M}$ is a retract of $F$ so the top-horizontal isomorphism in (\ref{everythinginsniceandfree}) descends to an isomorphism 
\[\pi_\ast^e\colon[\widetilde{M},N]_R^\gl\overset{\cong}{\longrightarrow} \hom_{\lmod_{\pi_\ast^eR}}(\pi_\ast^e \widetilde{M}, \pi_\ast^e N)\cong \hom_{\lmod_{\pi_\ast^eR}}(\pi_\ast^e M, \pi_\ast^e N).\]
Setting $N=M$ we then lift the isomorphism $\pi_\ast^e \widetilde{M}\cong \pi_\ast^e M$ to a map $\phi\in \hom_{\ho^\gl(\lmod_R)}(\widetilde{M},M)$. As both $M$ and $\widetilde{M}$ are globally flat and $\pi_\ast^e \phi$ is an isomorphism by construction, Remark~\ref{easyconsequence} says $\phi$ is an isomorphism inside $\ho^\gl(\lmod_R)$. The following commutative diagram of abelian groups then finishes our proof,
\[\begin{tikzcd}
{\left[M,N\right]_R^\gl}\ar[rr, "\pi_\ast^e"]\ar[d, "\phi^\ast"]\ar[d, "\cong", swap]	&&	{\hom_{\lmod_{\pi_\ast^e R}}(\pi_\ast^e M,\pi_\ast^e N)}\ar[d, "(\pi_\ast^e \phi)^\ast=\phi_g^\ast"]\ar[d, "\cong", swap]	\\
{[\widetilde{M},N]_R^\gl}\ar[rr, "\pi_\ast^e"]\ar[rr, "\cong", swap]		&&	{\hom_{\pi_\ast^e \lmod_R}(\pi_\ast^e \widetilde{M},\pi_\ast^e N)}	
\end{tikzcd}.\]
\end{proof}

A consequence of Proposition~\ref{proptwo} is the following strengthening of Proposition~\ref{extensiontoflatness}.

\begin{cor}\label{proponeagain}
Let $R$ be a global ring spectrum and $M_\ast$ a projective left $\pi_\ast^e R$-module. Then there exists a globally flat left $R$-module $M$ with $\pi_\ast^e M\cong M_\ast$, \emph{unique up to global equivalence}.
\end{cor}

\begin{proof}
Proposition~\ref{extensiontoflatness} gives us existence. Let $(M,\phi_f)$ and $(M', \phi_{f'})$ be two globally flat $R$-modules with isomorphisms $\phi_f\colon \pi_\ast^e M\cong M_\ast$ and $\phi_{f'}\colon \pi_\ast^e M'\cong M_\ast$, then lifting the isomorphism $\phi^{-1}_{f'}\circ \phi_f\colon \pi_\ast^e M\cong \pi_\ast^e M'$ using Proposition~\ref{proptwo} gives an isomorphism $\phi\colon M\rightarrow M'$ in $\ho^\gl(\lmod_R)$ by Remark~\ref{easyconsequence}.
\end{proof}

It was mentioned in \textsection\ref{homtopabsdo} that left induced left $R$-modules were hard to work with, in particular, their homotopy groups hard to calculate. The following theorem shows that left induced left $R$-modules are globally flat in special cases.

\begin{theorem}\label{superusefullater}
Let $R$ be a global ring spectrum, and $M$ a left $R$-module such that $\pi_\ast^e M$ is a projective left $\pi_\ast^e R$-module. Then $M$ is globally flat if and only if $M$ is left induced.
\end{theorem}

Let us use the same notation as in the proof of Corollary~\ref{somethingsimplethanksref}.

\begin{proof}
Without loss of generality we can take $M$ to be a $\gl$-fibrant $R$-module. By Corollary~\ref{somethingsimplethanksref}, it is necessary and sufficient to show an $e$-cofibrant replacement $c\colon M_c\to M$ is a global equivalence if and only if $M$ is globally flat.\\

As $\pi_\ast^e M$ is projective over $\pi_\ast^e R$, we use the proof of Proposition~\ref{extensiontoflatness}, with respect to the \emph{trivial} global family, to obtain a left $R$-module $\widetilde{M}=\hocolim_i^e M_i$, which is a sequential $e$-homotopy colimit of wedges of suspensions of $R$, such that there exists an isomorphism of left $\pi_\ast^e R$-modules $\phi_\ast\colon \pi_\ast^e \widetilde{M}\cong \pi_\ast^e M$. Using Proposition~\ref{proptwo}, again with respect to the trivial global family, we see the natural map
\[[\widetilde{M}, M_c]_R^e\xrightarrow{\pi_\ast^e}\hom_{\lmod_{\pi_\ast^e R}}(\pi_\ast^e \widetilde{M}, \pi_\ast^e M_c)\]
is an isomorphism of abelian groups, leading us to recognise $\phi_\ast$ by a morphism $\phi\colon \widetilde{M}\to M_c$ in $\ho^e(\lmod_R)$. As $\widetilde{M}$ is $e$-cofibrant and $M_c$ is $e$-fibrant, we can take $\phi$ to be a strict map in $\lmod_R^e$. As $\phi$ is an $e$-equivalence between $e$-cofibrant $R$-modules, then by Remark~\ref{cofibrantobjectsgivege} we see $\phi$ is a $\gl$-equivalence. Moreover, by Proposition~\ref{naturalnessandstuff} and the fact that sequential $e$-homotopy colimits and sequential $\gl$-homotopy colimits can be modelled by mapping telescopes, we see $\widetilde{M}$ is globally flat. The following naturality diagram of graded global functors shows $c$ is a gl-weak equivalence if and only if $M$ is globally flat,
\[\begin{tikzcd}
{\epi_\ast R\underset{\pi_\ast^e R}{\otimes} \pi_\ast^e {\widetilde{M}}}\ar[rr, "{\id\otimes (\pi_\ast^e \phi)}"]\ar[rr, "\cong", swap]\ar[d, "\cong", swap]\ar[d, "{\Lambda_{\widetilde{M}}}"]	&&	{\epi_\ast R\underset{\pi_\ast^e R}{\otimes} \pi_\ast^e M_c}\ar[rr, "{\id\otimes \pi_\ast^e c}"]\ar[rr, "\cong", swap]\ar[d, "{\Lambda_{M_c}}"]	&&	{\epi_\ast R\underset{\pi_\ast^e R}{\otimes}\pi_\ast^e M}\ar[d, "{\Lambda_{M}}"]	\\
{\epi_\ast \widetilde{M}}\ar[rr, "{\epi_\ast \phi}"]\ar[rr, "\cong", swap]	&&	{\epi_\ast M_c}\ar[rr, "{\epi_\ast c}"]	&&	{\epi_\ast M}
\end{tikzcd}.\]
\end{proof}

\begin{prop}\label{propthree}
Let $R$ be a global ring spectrum, $N$ a globally flat left $R$-module such that $\pi_\ast^e N$ is a projective left $\pi_\ast^e R$-module. Then for any right $R$-module $M$ there is an isomorphism of $\epi_\ast R$-modules
\[\epi_\ast M \underset{\pi_\ast^e R}{\otimes} \pi_\ast^e N\overset{\cong}{\longrightarrow} \epi_\ast (M\smplr N).\]
\end{prop}

\begin{proof}
The canonical map of this proposition is defined for each $G$ in $\F$ as
\[\Theta_{M,N}^G\colon \pi_\ast^G M\underset{\pi_\ast^e R}{\otimes} \pi_\ast^e N\longrightarrow \pi_\ast^G(M\smplr N),\qquad m\otimes n\longmapsto m{\times} p^\ast_G(n).\]
Above the operation $-\times -$ is the derived $R$-relative box product pairing, which is defined as follows: first one takes cofibrant replacements of $M$ and $N$, say $M_c$ and $N_c$, and then one considers the composition
\begin{equation}\label{bilinearmapping}\pi_\ast^G M_c \times \pi_\ast^G N_c\to \pi_\ast^G \left(M_c\otimes N_c\right)\to \pi_\ast^G\left(M_c\otimes_R N_c\right),\end{equation}
where the first morphism is the absolute box product pairing of \cite[Construction 3.5.12]{s}, and the second morphism is induced by postcomposed with the canonical map $M_c\otimes N_c\rightarrow M_c\otimes_R N_c$. This postcomposition and \cite[Theorem 3.5.14]{s} imply that (\ref{bilinearmapping}) is bilinear over $\pi_\ast^G R$, giving us the desired $\pi_\ast^GR$-linear derived $R$-relative box product
\[\pi_\ast^G M\underset{\pi_\ast^G R}{\otimes} \pi_\ast^GN\cong \pi_\ast^G M_c\underset{\pi_\ast^G R}{\otimes} \pi_\ast^G N_c\to \pi_\ast^G\left( M_c\otimes_R N_c\right)\cong \pi_\ast^G(M\smplr N).\]
The maps $\Theta_{M,N}^G$ have similar properties to the $\Lambda^G_M$ maps from Definition~\ref{globallyflat}, which is not remarkable as $\Lambda_M^G=\Theta_{R,M}^G$. The fact $\Theta^G_{M,N}$ is natural in the right $R$-module variable $M$ follows from the bifunctorality of $-\smpl_R-$, and the fact these maps are natural in $N$ and $G$ follow from the same reasoning of Proposition~\ref{naturalnessandstuff}. We now have a map of graded global functors
\[\Theta_{M,N}\colon \pi_\ast^G M\underset{\pi_\ast^e R}{\otimes} \pi_\ast^e N\longrightarrow \pi_\ast^G(M\smplr N).\]
Writing $\lmod_R^{\Theta}$ for the full subcategory of $\lmod_R^\gl$ consisting of left $R$-modules $N$ such that $\Theta_{M,N}$ is an isomorphism for all right $R$-modules $M$. One observes $R$ is in $\lmod_R^\Theta$ and $\lmod_R^\Theta$ is closed under arbitrary suspensions, wedges, and filtered homotopy colimits, using similar reasoning to Proposition~\ref{naturalnessandstuff} and the fact that $-\smpl_R-$ commutes these constructions as $\ho^\gl(\lmod_R)$ is a closed symmetric monoidal category. As $N$ is globally flat and $\pi_\ast^e N$ is a projective $\pi_\ast^e R$-module, Corollary~\ref{proponeagain} says $N$ is globally equivalent to an explicit model given in the proof of Proposition~\ref{extensiontoflatness}, i.e., as a sequential homotopy colimit of wedges of shifts of $R$, so such an $N$ is in $\lmod_R^\Theta$, which finishes our proof.
\end{proof}

\begin{remark}[Tor and Ext spectral sequences]\label{spectralsequences}
Propositions~\ref{proptwo} and~\ref{propthree} resemble degenerate cases of a global Ext and global Tor spectral sequences respectively (similar to those found in \cite[Section IV.4]{ekmm}). In fact these two statements would need to be used to construct such spectral sequences. This is done in \cite[Section 2.3]{me}, although the only practical application (according to the author) seems to be a weakening of the hypothesis of projectivity in Proposition~\ref{propthree} to a flatness hypothesis.
\end{remark}

The following is a generalisation of Proposition~\ref{extensiontoflatness} along the lines of \cite[Theorem 6]{wolbert}.

\begin{prop}\label{wolbertgenerlaiseation}
Let $R$ be a global ring spectrum and $M_\ast$ a left $\pi_\ast^e R$-module of projective dimension at most two such that for all $G$ in $\F$, the groups
\[\tor_1^{\pi_\ast^e R}(\pi_\ast^G R, M_\ast),\qquad \tor_2^{\pi_\ast^e R}(\pi_\ast^G R, M_\ast)\]
vanish. Then there exists a globally flat left $R$-module $M$ with $\pi_\ast^eM\cong M_\ast$.
\end{prop}

The necessity of the ``Tor-condition'' above will be clear in the proof, and in particular holds if $M_\ast$ or $\pi_\ast^G R$ is flat over $\pi_\ast^e R$.

\begin{proof}
This proof follows along the same lines as \cite[Theorem 6]{wolbert}. First we deal with the projective dimension one case, so let
\[0\longrightarrow P^1_\ast\overset{f}{\longrightarrow} P^0_\ast\longrightarrow M_\ast\longrightarrow 0\]
be a projective resolution of the left $\pi_\ast^e R$-module $M_\ast$. By Proposition~\ref{proptwo} and Corollary~\ref{proponeagain} we have globally flat left $R$-modules $P^1$ and $P^0$, and a map of $R$-modules $g\colon P^1\rightarrow P^0$ realising the first map in the projective resolution above. Define $M$ as the cofibre of $g$. To see $M$ is globally flat over $R$, consider the following diagram of abelian groups with exact rows, for each $G$ in $\F$
\begin{equation}\label{nothinglikeglobalflatness}\begin{tikzcd}
{\tor^{\pi_\ast^e R}_1(\pi_\ast^G R, M_\ast)}\ar[r]	&	{\pi_\ast^G R \underset{\pi_\ast^e R}{\otimes} P_\ast^1}\ar[d, "{\Lambda_{P^1}^G}"]\ar[d, "\cong", swap]\ar[r]	&	{\pi_\ast^G R \underset{\pi_\ast^e R}{\otimes} P_\ast^0}\ar[d, "{\Lambda_{P^0}^G}"]\ar[d, "\cong", swap]\ar[r]	&	{\pi_\ast^G R \underset{\pi_\ast^e R}{\otimes} M_\ast}\ar[d, "{\Lambda_{M}^G}"]\ar[r]	&	{0}	\\
{\cdots}\ar[r]	&	{\pi_\ast^G P^1}\ar[r]		&	{\pi_\ast^G P^0}\ar[r]		&	{\pi_\ast^G M}\ar[r]		&	{\cdots}
\end{tikzcd}.\end{equation}
By assumption, the $\tor_1$-group above vanishes, hence $g$ induces an injection on all global homotopy groups, from which we immediately obtain, for each $G$ inside $\F$, the short exact sequence of left $\pi_\ast^G R$-modules
\[0\longrightarrow \pi_\ast^G P^1\longrightarrow \pi_\ast^G P^0\longrightarrow \pi_\ast^G M\longrightarrow 0.\]
By (\ref{nothinglikeglobalflatness}) and the five lemma we see $M$ is globally flat. Setting $G=e$, we also obtain $\pi_\ast^e M\cong M_\ast$.\\

Suppose $M_\ast$ now has projective dimension two, or equivalently that we have two exact sequences
\begin{equation}\label{apearofsess} 0\longrightarrow P^2_\ast\longrightarrow P^1_\ast\longrightarrow Q_\ast\longrightarrow 0,\qquad\qquad 0\longrightarrow Q_\ast\overset{f}{\longrightarrow} P_\ast^0\longrightarrow M_\ast\longrightarrow 0,\end{equation}
where each $P_\ast^i$ is a projective left $\pi_\ast^e R$-module. Notice that the second short exact sequence above implies
\[\tor_1^{\pi_\ast^e R}(\pi_\ast^G R, Q_\ast)\cong \tor_2^{\pi_\ast^e R}(\pi_\ast^G R, M_\ast)=0.\]
Using the projective dimension one case above, we can realise the first sequence of (\ref{apearofsess}) by a cofibre sequence of globally flat left $R$-module spectra,
\begin{equation}\label{littlecofibresequence}P^2\longrightarrow P^1\longrightarrow Q.\end{equation}
We can also use Corollary~\ref{proponeagain} to obtain a globally flat left $R$-module $P^0$ recognising $P_\ast^0$. Consider the commutative diagram of abelian groups
\[\begin{tikzcd}
	&	{[Q, P^0]_R^\gl}\ar[d]\ar[r]	&	{[P^1, P^0]_R^\gl}\ar[d]\ar[r]	&	{[P^2, P^0]_R^\gl}\ar[d]	\\
{0}\ar[r]	&	{\hom_{\lmod_{\pi_\ast^e R}}(Q_\ast, P^0_\ast)}\ar[r]	&	{\hom_{\lmod_{\pi_\ast^e R}}(P^1_\ast, P^0_\ast)}\ar[r]	&	{\hom_{\lmod_{\pi_\ast^e R}}(P^2_\ast, P^0_\ast)}
\end{tikzcd}.\]
The top row exact by the cofibre sequence (\ref{littlecofibresequence}) and the bottom row by applying $\hom(-, P_\ast^0)$ to the first short exact sequence of (\ref{apearofsess}). Using Proposition~\ref{proptwo} for $P^1$ and $P^2$ mapping into $P^0$, we see the middle and right vertical maps are isomorphisms. A diagram chase then shows there is a map of left $R$-modules $g\colon Q\to P^0$ recognising $f$. Taking a fibrant replacement $P^0_f$ of $P^0$ in $\lmod_R^\gl$ ($Q$ is already cofibrant), we realise $g\colon Q\to P^0_f$ in $\lmod_R^\gl$ and define $M$ as the cofibre of this map. To see $M$ is globally flat, we use the same argument as in the projective dimension 1 case.
\end{proof}


\section{Realising algebra with homotopy global ring spectra}\label{sectionhomotopycommutative}
In this section we obtain our first realisation result with multiplicative structure.

\begin{theorem}\label{theoremone}
Let $R$ be an ultra-commutative ring spectrum and $\eta_\ast\colon \pi_\ast^e R\rightarrow S_\ast$ a map of graded commutative rings witnessing $S_\ast$ as a projective $\pi_\ast^e R$-module. Then there exists a globally flat homotopy commutative global $R$-algebra $S$ with $\pi_\ast^e S\cong S_\ast$, such that for all homotopy commutative $R$-algebras $T$ and all maps of $\pi_\ast^e R$-algebras $\psi_\ast\colon S_\ast\rightarrow \pi_\ast^e T$, there exists a unique map $\psi\colon S\rightarrow T$ of homotopy commutative global $R$-algebras such that $\eta_T= \psi\circ \eta_S$ inside $\ho^\gl(\mod_R)$. 
\end{theorem}

In particular, $S$ is the initial globally flat homotopy commutative global $R$-algebra recognising $S_\ast$ inside the homotopy category $\ho^\gl(\lmod_R)$.

\begin{remark}\label{generalisation}
The above theorem generalises to the case when $R$ is a cofibrant $\E_\infty$-global ring spectrum, as the only fact we need for the following proof is that the homotopy category $\ho^\gl(\mod_R)$ has a monoidal structure inherited from the (derived) smash product over $R$. If $R$ is a cofibrant $\E_\infty$-global ring spectrum, then write $\O$ for an $\E_\infty$-operad and by \cite[Proposition 2.3]{togetmodulesovereinfty} the \emph{enveloping algebra} $\mathrm{Env_\O(R)}$ (see \cite[Definition 1.11]{togetmodulesovereinfty}), is a well-pointed monoid in $\spec^\gl$. By \cite[Theorem 1.10]{togetmodulesovereinfty} the category of global $R$-modules, defined in the operadic sense (see \cite[Definition 1.1]{togetmodulesovereinfty}), is equivalent to the category of $\mathrm{Env_\O(R)}$-modules $\mod_{\mathrm{Env_\O(R)}}$. By \cite[Proposition 2.7(a)]{togetmodulesovereinfty} we see this category $\mod_{\mathrm{Env_\O(R)}}$ comes with a left induced model structure from $\spec^\gl$, which moreover has the expected monoidal structure. The monoidal structure on $\mod_{\mathrm{Env_\O(R)}}$ then induces the desired monoidal structure on $\ho^\gl(\mod_R)$.
\end{remark}

Recall that if $M$ is an $R$-module spectrum, then $M^{\otimes n}$ refers to the $n$-fold smash product of $M$ \emph{over $R$}, and similarly, for a $\pi_\ast^e R$-module $M_\ast$ the iterated tensor product is always over $\pi_\ast^e R$.

\begin{proof}[Proof of Theorem~\ref{theoremone}]
This proof follows along the same lines as \cite[Theorem 2.1.1]{br}. We obtain existence of $S$ by Corollary~\ref{proponeagain}, which gives us a globally flat $R$-module $S$ with $\pi_\ast^e S\cong S_\ast$. Proposition~\ref{propthree} iteratively calculates
\[\pi_\ast^G S^{\smp n}\cong \pi_\ast^G R\underset{\pi_\ast^e R}{\otimes} S_\ast^{\otimes n}\]
for all $G$ in $\F$ and Proposition~\ref{proptwo} gives us the first isomorphism
\[\Phi\colon \hom_{\ho^\gl(\mod_R)}(S^{\smp n}, N)\cong \hom_{\mod_{\pi_\ast^e R}}(\pi_\ast^e S^{\smp n}, \pi_\ast^e N)\cong\hom_{\mod_{\pi_\ast^e R}}(S_\ast^{\otimes n}, \pi_\ast^e N),\]
for all $R$-modules $N$. Setting $N=S$, we transport the unit and multiplication map of the $\pi_\ast^e R$-module $S_\ast$ along $\Phi$ for $n=0$ and $2$, respectively, to obtain unit and multiplication maps on $S$ inside $\ho^\gl(\mod_R)$. As $\mod_R$ is a monoidal model category and $S$ is bifibrant, $S^{\smp n}$ is also cofibrant by the pushout product axiom (see \cite[Definition 3.1]{ss}), and these unit and multiplication maps can be realised by strict maps of $R$-modules $\eta\colon R\rightarrow S$ and $\mu\colon S^{\smp 2}\rightarrow S$. The unitality, associativity, and commutativity of these maps in $\ho^\gl(\mod_R)$ come from $\Phi$ for $n=1$, $3$, and $2$, respectively, again setting $N=S$.\\

To show the existence and uniqueness of $\psi$, let $T$ be a homotopy commutative $R$-algebra and $\phi_\ast\colon S_\ast\rightarrow \pi_\ast^e T$ a map of $\pi_\ast^e R$-algebras. Recall the set of morphisms of commutative monoids in a symmetric monoidal category $(\C, \otimes, \1)$ can be written as the equaliser
\[\hom_{\calg(\C)}(A,B)\longrightarrow \hom_\C(A,B) \rightrightarrows \hom_\C(A^{\otimes 2}, B)\times \hom_\C(\1, B),\]
where the parallel maps send $f\mapsto (f\circ \mu_A, f\circ \eta_A)$ and $f\mapsto (\mu_B\circ(f\otimes f), \eta_B)$, and $\calg(\C)$ denotes the category of commutative algebra objects of $\C$. Applying this to the symmetric monoidal categories $(\ho^\gl(\mod_R), \smpl_R, R)$ and $(\mod_{\pi_\ast^e R}, \otimes_{\pi_\ast^e R}, \pi_\ast^e R)$ with $A=S$ and $B=T$, and using $\Phi$ with $N=T$ we obtain the natural bijection
\begin{equation}\label{morehelpfulguy}\hom_{\calg(\ho^\gl(\mod_R))}(S,T)\overset{\cong}{\longrightarrow} \hom_{\calg_{\pi_\ast^e R}}(S_\ast, \pi_\ast^e T).\end{equation}
This allows us to lift $\psi_\ast\colon S_\ast\rightarrow \pi_\ast^e T$ to a unique map $\psi\colon S\rightarrow T$ in $\calg (\ho^\gl(\mod_R))$.
\end{proof}


\section{Realising algebra with $\E_\infty$-global ring spectra}\label{operads}
Using nonequivariant obstruction theory, we can place an $\E_\infty$-structure on the $S$ in Theorem~\ref{theoremone}, given some more conditions on $\eta_\ast\colon \pi_\ast^e R\rightarrow S_\ast$. To access this nonequivariant obstruction theory, we need some statements about endomorphism operads. Recall the model structure on the category of topological operads from \cite[Example 3.3.2]{bm}, where weak equivalences and fibrations are given level-wise.

\begin{lemma}\label{babylemma}
Let $\M$ be a topological monoidal model category. If $f\colon X\rightarrow Y$ is an acyclic fibration between bifibrant objects, then there is a zigzag of weak equivalences of topological operads
\[\opend_\M(X)\simeq\opend_\M(Y).\]
\end{lemma}

\begin{proof}
Define a topological operad $\opend_\M(f)$ at level $n$ by the following pullback diagram of spaces,
\begin{equation}\label{niceoperaddiagram}\begin{tikzcd}
{\opend_\M(f)_n}\ar[rr, "\pi_Y"]\ar[d, "\pi_X"]	&&	{\Map_\M(Y^{\otimes n},Y)=\opend_\M(Y)_n}\ar[d, "{(f^{\otimes n})^\ast=\Map_\M(f^{\otimes n}, Y)}"]	\\
{\Map_\M(X^{\otimes n},X)=\opend_\M(X)_n}\ar[rr, "{f_\ast=\Map_\M(X^{\otimes n}, f)}"]	&&	{\Map_\M(X^{\otimes n}, Y)}
\end{tikzcd}.\end{equation}
The composition operation on $\opend_\M(f)$ is the product of the composition operations on $\opend_\M(X)$ and $\opend_\M(Z)$, and in this way $\pi_X$ and $\pi_Y$ induce maps of topological operads. To be a little more precise, given two nonnegative integers $m$ and $n$, the composition operation
\[\opend_\M(f)_m\times \opend_\M(f)_n\times \opend_\M(f)_2\to \opend_\M(f)_{m+n}\]
is explicitly given by the assignment
\[\left((X^{\otimes m}\xrightarrow{g_m} X, Y^{\otimes m}\xrightarrow{h_m} Y), (X^{\otimes n}\xrightarrow{g_n} X, Y^{\otimes n}\xrightarrow{h_n} Y), (X^{\otimes 2}\xrightarrow{g} X, Y^{\otimes 2}\xrightarrow{h} Y)\right)\]
\[\longmapsto \left(X^{\otimes (m+n)}\xrightarrow{g_m\otimes g_n} X^{\otimes 2}\xrightarrow{g} X, Y^{\otimes(m+n)}\xrightarrow{h_1\otimes h_2}Y^{\otimes 2}\xrightarrow{h} Y \right).\]
This composition operation generalises to arbitrary $n$-tuples of nonnegative integers in the obvious way. From this definition, it is clear that $\pi_X$ and $\pi_Y$ commute with the various composition operations on $\opend_\M(X)$, $\opend_\M(Y)$, and $\opend_\M(f)$, inducing morphisms of topological operads.\\

As $f$ is an acyclic fibration, $X^{\otimes n}$ is cofibrant, and $X$ and $Y$ are fibrant, we see $f_\ast$ is also an acyclic fibration of spaces. Similarly, as $f^{\otimes n}$ is a weak equivalence, $X^{\otimes n}$ and $Y^{\otimes n}$ are cofibrant, and $Y$ is fibrant, we see $(f^{\otimes n})^\ast$ is a weak homotopy equivalence of spaces. We conclude $\pi_X$ is a weak equivalence as the category of topological spaces is (right) proper, and $\pi_Y$ is also a weak equivalence (an acyclic fibration even) as a base change of an acyclic fibration. As $\pi_X$ and $\pi_Y$ assemble to form maps of topological operads, the above argument witnesses these assembled maps as weak equivalences of topological operads. Hence we obtain a zigzag of weak equivalences
\[\opend_\M(X) \overset{\simeq}{\longleftarrow} \opend_\M(f) \overset{\simeq}{\longrightarrow} \opend_\M(Y).\qedhere\]
\end{proof}

Let $\M_1$ and $\M_2$ be two model categories with the same underlying category. Let $i$-$\mathcal{P}$ be the adjective referring to the model categorical property $\mathcal{P}$ inside $\M_i$, for $i=1,2$.

\begin{theorem}\label{bifibrantendotheorem}
Let $\M_1$ and $\M_2$ be topological monoidal model categories with the same underlying symmetric monoidal category $\M$ such that the $1$-weak equivalences are contained in the $2$-weak equivalences and the $1$-fibrations are contained in the $2$-fibrations. If $f\colon X\rightarrow Y$ is a 1-weak equivalence, where $X$ is 2-bifibrant and $Y$ 1-bifibrant, then there is a zigzag of weak equivalences between the topological endomorphism operads
\[\opend_{\M_2}(X)\simeq\opend_{\M_1}(Y).\]
In particular, if $\M_1=\M_2$ as model categories, then a weak equivalence $f\colon X\rightarrow Y$ between bifibrant objects induces a zigzag of weak equivalences between endomorphism operads.
\end{theorem}

\begin{proof}
First we factorise $f$ as a 1-acyclic cofibration followed by a 1-acyclic fibration
\[ X\overset{i}{\longrightarrow} Z \overset{p}{\longrightarrow} Y .\]
Notice the 1-acyclic fibrations are contained in the 2-acyclic fibrations, so by lifting properties we see 2-cofibrations are contained inside 1-cofibrations. In particular, 2-cofibrant objects are 1-cofibrant. We then see that $Z$ is 1-bifibrant as $X$ is 1-cofibrant and $Y$ is 1-fibrant. We now define a topological operad $\opend(i)$ at level $n$ by the following pullback diagram of spaces,
\begin{equation}\label{niceoperaddiagram}\begin{tikzcd}
{\opend(i)_n}\ar[rrr, "\pi_Z"]\ar[d, "\pi_X"]	&&&	{\Map_{\M_1}(Z^{\otimes n},Z)}\ar[d, "{(i^{\otimes n})^\ast=\Map_{\M_1}(Z^{\otimes n}, Z)}"]	\\
{\Map_{\M_2}(X^{\otimes n},X)}\ar[rrr, "{i_\ast=\Map_{\M_2}(X^{\otimes n}, i)}"]	&&&	{\Map_{\M_2}(X^{\otimes n}, Z)=\Map_{\M_1}(X^{\otimes n},Z)}
\end{tikzcd}.\end{equation}
Similar to the proof of Lemma~\ref{babylemma}, the composition operation on $\opend(i)$ is the product of that on $\opend_{\M_2}(X)$ and $\opend_{\M_1}(Z)$ such that $\pi_X$ and $\pi_Z$ both induce morphisms of topological operads. The product map $i^{\otimes n}$ is a 1-acyclic cofibration by the pushout product axiom, and $Z$ is $1$-fibrant, so $(i^{\otimes n})^\ast$ is an acyclic fibration of spaces. Similarly, $i_\ast$ is weak homotopy equivalence of spaces, as $X^{\otimes n}$ is 2-cofibrant, $X$ and $Z$ are 2-fibrant, and $i$ is a $2$-weak equivalence. Similar to Lemma~\ref{babylemma}, we see $\pi_X$ and $\pi_Z$ are both weak homotopy equivalences of spaces. This gives us the zigzag of weak equivalences of topological operads
\[\opend_{\M_2}(X) \overset{\simeq}{\longleftarrow} \opend(i) \overset{\simeq}{\longrightarrow} \opend_{\M_1}(Z).\]
Using Lemma~\ref{babylemma} with respect to $p$ we obtain the following zigzag of weak equivalences of topological operads
\[\opend_{\M_1}(Z) \overset{\simeq}{\longleftarrow} \opend_{\M_1}(p) \overset{\simeq}{\longrightarrow} \opend_{\M_1}(Y).\]
Combining the two zigzags above, we obtain the desired result.
\end{proof}

Setting $\M_1=\M_2$ in Theorem~\ref{bifibrantendotheorem} one obtains a generalisation of Lemma~\ref{babylemma} to the case when $f$ is simply a weak equivalence between bifibrant objects. There is also a dual statement to Theorem~\ref{bifibrantendotheorem}.

\begin{cor}
Let $\M_1$ and $\M_2$ be topological monoidal model categories, with the same underlying monoidal category $\M$, such that the $1$-weak equivalences are contained in the $2$-weak equivalences, and the $1$-cofibrations are contained in the $2$-cofibrations. If $f\colon X\rightarrow Y$ is a 1-weak equivalence, where $X$ is 1-bifibrant and $Y$ 2-bifibrant, then there is a zigzag of weak equivalences between the topological endomorphism operads
\[\opend_{\M_1}(X)\simeq\opend_{\M_2}(Y).\]
\end{cor}

\begin{proof}
First we factorise $f$ as a 1-acyclic cofibration followed by a 1-acyclic fibration. The result follows from the ``in particular'' statement of Theorem~\ref{bifibrantendotheorem} for the 1-acyclic cofibration. The rest of the proof is dual to the proof of Theorem~\ref{bifibrantendotheorem}.
\end{proof}

Before we prove the main result of this section, let us recall \cite[Definition 3.1]{r}, which we relativise over a base $\E_\infty$-ring $A$. Notice this is a purely nonequivariant condition.

\begin{mydef}\label{perfectuniversalblabla}
Let $A$ be an $\E_\infty$-ring and $B$ a homotopy commutative $A$-algebra. Then we say the $A$-algebra $B$ satisfies the \emph{perfect universal coefficient formula} if the following two conditions hold:
\begin{enumerate}
\item The graded ring $\pi_\ast (B\smpl_A B)$ is flat over $\pi_\ast B$.
\item For every $n\geq 1$, the natural map
\[[B^{\otimes n}, B]_A\xrightarrow{\pi_\ast}\hom_{\mod_{\pi_\ast A}}(\pi_\ast(B^{\otimes n}), \pi_\ast B)\]
is an isomorphism, where the (derived) smash products above are taken relative to $A$.
\end{enumerate}
\end{mydef}

We can now state the main theorem of this section. Recall the definition of a \'{e}tale morphism of (graded) commutative rings from \cite[\href{https://stacks.math.columbia.edu/tag/00U0}{Tag 00U0}]{stacks}.

\begin{theorem}\label{theoremthree}
Let $R$ be an ultra-commutative ring spectrum and $S$ a homotopy commutative global $R$-algebra which is left induced as an $R$-module and which satisfies the perfect universal coefficent formula as a nonequivariant homotopy commutative $R$-algebra. Suppose that either $\pi_\ast^e R\to \pi_\ast^e S$ is an \'{e}tale morphism of graded commutative rings or that it is a localisation.\footnote{The sentence ``Suppose that\ldots'' can be replaced by any sentence which implies that the $\Gamma$-cotangent complex $\mathcal{K}(B/A; M)$ of \cite[Section 3.2]{rw} is contractible, and the conclusion of the theorem will remain valid.} Then $S$ has an $\E_\infty$-global $R$-algebra structure, unique up to contractible choice in $\mod_R^\gl$, and the natural map
\begin{equation}\label{mappingpsacesstuff}\Map_{\calg_R^\gl}(S,S)\xrightarrow{\pi_\ast^e} \hom_{\calg_{\pi_\ast^e R}}(\pi_\ast^e S, \pi_\ast^e S)\end{equation}
is a weak equivalence of spaces, where the codomain is discrete.
\end{theorem}

Unique up to contractible choice means a certain moduli space (\`{a} la \cite{gh04}) is contractible.

\begin{remark}\label{generalisation2}
Just as in Remark~\ref{generalisation}, the above theorem generalises to the case when $R$ is simply a cofibrant $\E_\infty$-global ring spectrum. We suggest the interested reader follows the proof of Theorem~\ref{theoremthree} in the situation when $R$ is ultra-commutative, and then comes back to this remark. Indeed, suppose we are in the situation of Theorem~\ref{theoremthree} where $R$ is only assumed to be an $\E_\infty$-global ring spectrum. First, we take a cofibrant replacement of $R$ inside $\calg^\gl$, using the model structure of Definition~\ref{einf}. An $\E_\infty$-global $R$-algebra structure on $S$ is an $\E_\infty$-structure on $S$ in $\spec^\gl$ and a morphism $R\to S$ in $\calg^\gl$. By \cite[Lemma 1.7]{togetmodulesovereinfty}, we see this is equivalent to an $\O_R$-structure on $S$ in $\spec^\gl$, where $\O_R$ is the enveloping operad of the $\O$-algebra $R$, where $\O$ is an $\E_\infty$-operad; see \cite[Definition 1.5]{togetmodulesovereinfty}. One can then replace each occurance of $\O$ (resp. $\opend_R^-(-)$) in the whole of the proof of Theorem~\ref{theoremthree} below with $\O_R$ (resp. $\opend_\sph^-(-)$), using the facts that $\O$ is $\Sigma$-cofibrant and $R$ is cofibrant in $\calg^\gl$ in tandem with \cite[Proposition 2.3]{togetmodulesovereinfty} to see $\O_R$ is an admissible and $\Sigma$-cofibrant operad. Notice the second half of the proof (regarding the nonequivariant deformation theory) remains untouched by this change.
\end{remark}

The proof of Theorem~\ref{theoremthree} uses the vanishing of certain obstruction groups found in \cite{r}. The \'{e}tale case can be found in \cite{rw}, and the localisation case we do ourselves now.

\begin{lemma}\label{localisationkillsgammacotangntcomplex}
Let $A$ be a graded commutative ring (considered as an $\E_\infty$-dga with trivial differential) and $S\subseteq A$ be a multiplicative subset of $A$. Then for every graded $A[S^{-1}]$-module $M$, the $\Gamma$-cotangent complex $\mathcal{K}(A[S^{-1}]/A; M)$ is contractible in the derived category of $A[S^{-1}]$.
\end{lemma}

\begin{proof}
The proof follows the analogous argument for the usual cotangent complex of Quillen; see \cite[Proposition 5.1]{quillencotangentcomplex}. Writing $B=A[S^{-1}]$, then a simple fact about localisation is that $B\otimes_A C\simeq B\smpl_A C$ is quasi-isomorphic to $C$ for every $B$-complex $C$. This fact and the flat base change of \cite[Theorem 5.8(1)]{rw} give us the following chain of quasi-isomorphisms
\[\mathcal{K}(B/A;M)\simeq \mathcal{K}(B/A;M)\smpl_A B\simeq \mathcal{K}(B\smpl_A B/B;M)\simeq \mathcal{K}(B/B;M)\simeq 0,\]
where the last quasi-isomorphism comes from the definition \cite[Paragraph 3.2]{rw}.
\end{proof}

\begin{proof}[Proof of Theorem~\ref{theoremthree}]
Using the same notation as the proof of Corollary~\ref{somethingsimplethanksref}, we can without loss of generality take $R$ to be $\gl$-bifibrant, and the fact $S$ is left induced means an $e$-cofibrant replacement of $R$-modules $S_c\to S$ is a $\gl$-equivalence; see Corollary~\ref{somethingsimplethanksref}. As done in the proof of \cite[Proposition 2.2.3]{br}, we use a relativised version of \cite[Corollary 5.8]{r} (using our perfect universal coefficient formula assumption), and by either \cite[Theorem 5.8(3)]{rw} in the \'{e}tale case or Lemma~\ref{localisationkillsgammacotangntcomplex} in the localisation case, we obtain an $e$-$\E_\infty$-structure on $S_c$. In other words, we obtain an $\E_\infty$-structure on $S_c$ inside $\mod_R^e$, so a map of topological operads
\[\ga\colon \O\longrightarrow \opend_R^e(S_c),\]
where $\O$ is an $\E_\infty$-operad. We are now in the position to use Theorem~\ref{bifibrantendotheorem} with respect to $\M_1=\mod_R^\gl$, $\M_2=\mod_R^e$ (see Remark~\ref{modelstructurestwo}), and $f\colon S_c\longrightarrow S$, which gives us a zigzag of weak equivalences of topological operads
\begin{equation}\label{greatness}\opend_R^e(S_c)\simeq \opend_R^\gl(S).\end{equation}
In particular, we obtain a bijection of sets
\[[\O, \opend_R^e(S_c)]_\topop\cong[\O,\opend_R^\gl(S)]_\topop,\]
where $\topop$ denotes the category of topological operads with the model structure of \cite[Example 3.3.2]{bm}. Using the fact that $\O$ is cofibrant in $\topop$ and all objects are fibrant, we define our $\gl$-$\E_\infty$-structure on $S$ to be the image of $\ga$ under the above isomorphism.\\

To show this $\E_\infty$-structure is unique up to homotopy, observe the following chain of bijections of $\pi_0$ of (derived) mapping spaces of topological operads,
\begin{equation}\label{oneconnectedspace}\pi_0 \Map_{\topop}(\O, \opend_R^\gl(S))\cong \pi_0\Map_{\topop}(\O, \opend_R^e(S_c))= \ast.\end{equation}
The first isomorphism is induced by (\ref{greatness}), and the second by \cite[Theorem 5.8(3)]{rw} and either \cite[Corollary 5.8]{r} in the \'{e}tale case and Lemma~\ref{localisationkillsgammacotangntcomplex} in the localisation case. At this stage we use $\ga$ to view $S$ (resp. $S_c$) as an object of $\calg^\gl_R$ (resp. $\calg_R^e$). Let $S^\O$ be a cofibrant replacement of $S$ in $\calg_R^\gl$, and $S^\O_c$ for a cofibrant replacement of $S_c$ in $\calg_R^e$. Consider the composition
\begin{equation}\label{newreplacement}S^\O_c\xrightarrow{\simeq} S_c\xrightarrow{\simeq} S\xleftarrow{\simeq} S^\O.\end{equation}
The first map is an $e$-equivalence between cofibrant nonequivariant $\E_\infty$-$R$-algebras (hence cofibrant nonequivariant $R$-modules) and hence a $\gl$-equivalence by Remark~\ref{cofibrantobjectsgivege}, and the second map is a $\gl$-equivalence as $S$ is left induced (see Corollary~\ref{somethingsimplethanksref}), hence the composition (\ref{newreplacement}) is a global equivalence. By the usual arguments, $S^\O_c$ is bifibrant in $\calg_R^e$, hence cofibrant in $\calg_R^\gl$,  and $S^\O$ is bifibrant in $\calg_R^\gl$, hence (\ref{newreplacement}) can be realised by a single map in $\calg_R^\gl$. Considering these replacements now, we drop the superscript $\O$ from our notation.\\

To see the $\E_\infty$-$R$-algebra structure on $S$ is unique up to contractible choice, we need to define a moduli space, which we do following \cite[Section 5]{gh04}. Considering $\calg_R^\gl$ as a simplicial model category via the singular set functor, we let $\M_R^\gl(S)$ be the classifying space of the category $\mathcal{E}(S)$ of $\E_\infty$-global $R$-algebras $T$ which are isomorphic to $S$ inside the category of commutative algebra objects of $\ho^\gl(\mod_R)$ (the isomorphism is \textbf{not} part of the data) and with morphisms that are $\gl$-equivalences. By \cite{dwyerkannightmare} we see there are weak equivalences of spaces
\begin{equation}\label{dwyerkangavemthis}\M_R^\gl(S)\simeq \coprod_{[T]}\M(T)\simeq \coprod_{[T]} B\mathrm{Aut}(T),\end{equation}
where the coproduct is indexed by global equivalence classes of objects $T$ in $\mathcal{E}(S)$ (this is a set using the definitions of \cite{dwyerkannightmare}), $\M(T)$ is the classifying space of the subcategory of $\calg_R^\gl$ consisting of objects equivalent to a chosen (bifibrant) representative $T$ and $\gl$-equivalences, and $\mathrm{Aut}(T)$ is the monoid component of automorphisms of $T$ in $\calg_R^\gl$. Let us first notice that the space $\M_R^\gl(S)$ is nonempty and path-connected by (\ref{oneconnectedspace}), so $\M_R^\gl(S)$ is contractible if and only if $\Omega_\ga \M_R^\gl(S)$ is contractible. From (\ref{dwyerkangavemthis}) we see that
\[\Omega_\ga \M_R^\gl(S)\simeq \Omega B \mathrm{Aut}(S)\simeq \mathrm{Aut}(S)\subseteq \Map_{\calg^\gl_R}(S, S)\]
is the path component of $\Map_{\calg_R^\gl}(S, S)$ based at the identity. From the fibrancy conditions compiled above for $S$ and $S_c$, then the fact $S$ is left induced gives us a chain of weak equivalences of derived mapping spaces
\begin{equation}\label{mappingspacesarewelldefined}\Map_{\calg_R^\gl}(S,S)\simeq \Map_{\calg_R^\gl}(S_c, S)\simeq \Map_{\calg_R^e}(S_c, S)\simeq \Map_{\calg_R^e}(S_c, S_c)=\overline{\M},\end{equation}
the second weak equivalence in induced by the Quillen adjunction between $\calg_R^\gl$ and $\calg_R^e$ given by the identity. Our goal now is to show $\overline{\M}$ is discrete. Following the proof of \cite[Theorem 2.2.4]{br}, we use a relativised version of \cite[Theorem 4.5]{gh04} over $R$, with $X=Y=E=S$ (the fact $S$ satisfies the perfect universal coefficient formula as an $R$-algebra implies the Adams condition required in \cite[Definition 3.1]{gh04}). From this we obtain a second quadrant spectral sequence converging to the homotopy groups of $\overline{\M}$ based at the identity,
\[E_2^{s,t}\cong\left\{ \begin{array}{cc} \hom_{\calg_{\pi_\ast^e S}}(\pi_\ast^e(S\smplr S), \pi_\ast^e S) & (s,t)=(0,0) \\ \mathrm{Der}^s_{\pi_\ast^e S}(\pi_\ast^e(S\smplr S), (\pi_\ast^e S)[-t]) & t>0  \end{array}\right\} \Longrightarrow \pi_{t-s}\overline{\M},\]
where the homomorphism set is that of graded commutative $\pi_\ast^e S$-algebras, and $\mathrm{Der}^s_{A_\ast}(B_\ast, C_\ast[-t])$ is the $s$th derived functor of $A_\ast$-linear derivations into $C_{\ast+t}$ (see \cite[Section 4]{gh04}). Using the fact that $\mathcal{K}(\pi_\ast^e S/\pi_\ast^e R;M)$ is contractible for all $\pi_\ast^e S$-modules $M$, then the comparison results of \cite{comparisonresults} show that the above $E_2$-page is concentrated in filtration $t=0$, meaning it collapses on the $E_2$-page, and shows $\overline{\M}$ is weakly equivalent to a discrete space with 
\[\pi_0 \overline{\M}\cong \hom_{\calg_{\pi_\ast^e S}}(\pi_\ast^e(S\smplr S), \pi_\ast^e S)\cong \hom_{\calg_{\pi_\ast^e R}}(\pi_\ast^e S, \pi_\ast^e S).\]
In particular, we see that $\M_R^\gl(S)$ is contractible, hence the $\E_\infty$-global $R$-algebra structure on $S$ is unique up to contractible choice. Moreover, this argument and (\ref{mappingspacesarewelldefined}) show (\ref{mappingpsacesstuff}) is an isomorphism.
\end{proof}

Theorem~\ref{theoremthree} ties nicely into our continuing story about realising objects in global homotopy theory straight from algebraic information.

\begin{cor}\label{onlyusedforglaoisextensions}
Suppose we are in the situation of Theorem~\ref{theoremone}. If $\eta_\ast\colon \pi_\ast^e R\rightarrow S_\ast$ is in addition an \'{e}tale morphism of graded commutative rings, then the globally flat homotopy commutative $R$-algebra $S$ of Theorem~\ref{theoremone} realising $S_\ast$, has an $\E_\infty$-global $R$-algebra structure, unique up to contractible choice, and the natural map
\[\Map_{\calg_R^\gl}(S,S)\xrightarrow{\pi_\ast^e} \hom_{\calg_{\pi_\ast^e R}}(S_\ast, S_\ast)\]
is a weak equivalence of spaces, where the codomain is discrete.
\end{cor}

\begin{proof}
The proof of Theorem~\ref{theoremone} states that $S$ can be modelled by a bifibrant globally flat $R$-module, and Theorem~\ref{superusefullater} states that the projectivity of $\pi_\ast^e S$ over $\pi_\ast^eR$ implies $S$ is also left induced. Moreover, Propositions~\ref{proptwo} and~\ref{propthree} show that $S$ satisfies the perfect universal coefficient formula as an nonequivariant $R$-algebra (recall finite relative tensor products of projective modules are projective modules), placing us within the hypotheses of Theorem~\ref{theoremthree} above.
\end{proof}


\section{Realising algebra with $\G_\infty$-ring spectra}\label{ginfinity}
After Theorem~\ref{theoremthree}, one might have the following query:
\begin{center}	\emph{Why have we not placed an ultra-commutative structure on the $S$ from Theorem~\ref{theoremthree} despite the fact $R$ is an ultra-commutative ring spectrum?}	\end{center}
The answer is that the obstruction theory for ultra-commutative ring spectra akin to the $\E_\infty$-obstruction theory of \cite{gh04} and \cite{r} has not been developed yet. However this section aims to find a compromise.\\

The difference between ultra-commutative ring spectra and $\E_\infty$-global ring spectra that one can detect on their homotopy groups is the presence of power operations; see \cite[Definition 5.1.6, Theorem 5.1.11]{s}. The concept of a $\G_\infty$-structure on global spectra is discussed in \cite[Remark 5.1.16]{s} and in-depth in \cite{m}, and is the minimal structure on a global homotopy type to have power operations.  A $\G_\infty$-spectrum in global homotopy theory is analogous to an $H_\infty$-spectrum in classical homotopy theory; see \cite[Remark 5.1.14]{s}. Let us first define the spectra we need.

\begin{construction}
Let $G$ be a compact Lie group inside $\F$. We define the global spectra $\Sigma_+^\infty B_\gl G$ as $\Sigma_+^\infty \mathbf{L}_{G,V}$ where $V$ is any faithful $G$-representation; see \cite[Construction 1.1.27, Construction 4.1.7]{s}. By \cite[Proposition 1.1.26]{s} this is well-defined up to a preferred zigzag of global equivalences; see \cite[Definition 1.1.27]{s}. For any $m\geq 0$ we define $\Sigma_+^\infty E_\gl \Sigma_m$ to be the orthogonal spectrum $\Sigma_+^\infty \mathbf{L}_{\Sigma_m, \mathbb{R}^m}$, where $\mathbb{R}^m$ has the tautological $\Sigma_m$-action; see \cite[p.27, Construction 4.1.7]{s}. It follows from \cite[Proposition 1.1.26]{s} that $\Sigma_+^\infty E_\gl \Sigma_m$ is globally contractible. Notice that the $\Sigma_m$-coinvariants of $\Sigma_+^\infty E_\gl \Sigma_m$ are precisely $\Sigma_+^\infty B_\gl \Sigma_m$ by \cite[Definition 1.1.27]{s}, and $\Sigma_+^\infty B_\gl \Sigma_m$ has the nonequivariant homotopy type of $\Sigma_+^\infty B\Sigma_m$; see \cite[Remark 1.1.29]{s}.
\end{construction}

We will also need to recall the following general construction. 

\begin{construction}\label{grouprings}
Let $G$ be a finite group and $(\C,\otimes,\1)$ a closed symmetric monoidal category with finite coproducts. For a monoid $C$ of $\C$ we obtain a monoid $C[G]=\coprod_G C$, whose multiplication is defined through the multiplication on $C$ and the group $G$. In particular, if $\C=\spec^\gl$ and $R$ is a global ring spectrum, then $R[G]$ is a global ring spectrum with $\epi_\ast (R[G])\cong (\epi_\ast R)[G]$.
\end{construction}

Recall again, that iterated tensor products of $R$-modules are taken relative to $R$, for both module spectra and algebraic modules.

\begin{mydef}\label{ginfintyguys}
Let $R$ be an ultra-commutative ring spectrum and $M$ an $R$-module. For a fixed $1\leq n\leq \infty$, a \emph{$\G_n$-structure on $M$} is a series of maps in $\ho^\gl(\mod_R)$, for all $1\leq m\leq n$,
\[h_m\colon \mathbb{L}\P^m_R M\longrightarrow M, \qquad \mathbb{L}\P^m_R M = R\smp \left(\Sigma_+^\infty E_\gl \Sigma_m\underset{\Sigma_m}{\otimes} M^{\otimes m}\right)\cong \left(R\smp \Sigma_+^\infty E_\gl \Sigma_m\right)\underset{R[\Sigma_m]}{\otimes} M^{\otimes m},\]
such that for all integers $i,j, k,l$, with $i+j\leq n$ and $kl\leq n$, the following diagrams (from \cite[Proposition 1.12]{m}) commute in $\ho^\gl(\mod_R)$,
\[\begin{tikzcd}
{\mathbb{L}\P_R^iM\otimes_R \mathbb{L}\P_R^jM}\ar[r, "h_i\otimes h_j"]\ar[d]	&	{M\otimes_R M}\ar[r]	&	{\mathbb{L}\P^2 M}\ar[d, "h_2"]	\\
{\mathbb{L}\P_R^{i+j}M}\ar[rr, "h_{i+j}"]	&&	{M}
\end{tikzcd},\qquad
\begin{tikzcd}
{\mathbb{L}\P^k_R\mathbb{L}\P^l_R M}\ar[r, "{\P^k_Rh_l}"]\ar[d]	&	{\mathbb{L}\P^k_R M}\ar[d, "h_k"]	\\
{\mathbb{L}\P^{kl}_R M}\ar[r, "h_{kl}"]	&	{M}
\end{tikzcd}.\]
\end{mydef}

We justify the use of the notation $\mathbb{L}\P^m_R$ by \cite[Theorem 3.30]{m}, which states the definition of $\mathbb{L}\P^m_R$ above is a model for the left derived functor of the symmetric $R$-algebra functor $\P_R=\bigvee_{m\geq 0}(-)^{\otimes m}/\Sigma_m$.

\begin{theorem}\label{theoremfive}
Let $R$ be an ultra-commutative ring spectrum and $\eta_\ast\colon \pi_\ast^e R\rightarrow S_\ast$ a map of graded commutative rings which witnesses $S_\ast$ as a projective $\pi_\ast^e R$-module. If $S_0$ is a $\Z[1/n!]$-algebra for some $n\geq 1$ then there exists a globally flat homotopy commutative $R$-algebra $S$ such that $\pi_\ast^e S\cong S_\ast$, with a unique (up to global equivalence) $\G_n$-structure lifting the homotopy commutative multiplication on $S$. In particular, if $S_0$ is a $\Q$-algebra, then $S$ is a $\G_\infty$-$R$-algebra. 
\end{theorem}

\begin{remark}\label{generaliseaointthree}
Similar to Remark~\ref{generalisation}, the proof of Theorem~\ref{theoremfive} also holds in the more general case that $R$ is an $\E_\infty$-global ring spectrum. This might seem a little surprising, because the statement of the above theorem seems to imply that our $R$-module $S$ inherits its power operations from the ultra-commutative ring spectrum $R$, however, this is a red herring. Indeed, in the proof below it is clear that the $\G_n$-structure on $S$ (i.e.~the power operations) comes from the fact that $S_0$ is a $\Z[1/n!]$-algebra, \emph{not} the power operations on $R$.
\end{remark}

We will use a small lemma from homological algebra to obtain the above statement.

\begin{lemma}\label{littlehomologicalalgebra}
Let $R$ be a graded commutative ring, $M$ a graded $R$-module, and $m$ a positive integer. Suppose that each $M_n$, the submodule of $M$ concentrated in degree $n\in \Z$, is a $\Z[1/m!]$-module. If $M$ is a projective as a graded $R$-module, then $M^{\otimes m}$ is a projective left $R[\Sigma_m]$-module, and $(M^{\otimes m})_{\Sigma_m}$ is a projective $R$-module.
\end{lemma}

\begin{proof}
We will prove these facts in the opposite order. The tensor-hom adjunction shows inductively that if $M$ is a projective $R$-module then any tensor power of $M$ over $R$ is projective as an $R$-module. In general, if a finite group $H$ acts on an $R$-module $M$ by $R$-module homomorphisms, then as long as $M$ is a module over $\Z[1/|H|]$ the canonical map into the $H$-coinvariants $M\rightarrow M_H$ has a splitting
\[M_H\longrightarrow M,\qquad [x]\longmapsto \frac{1}{|H|}\sum_{h\in H} xh.\]
In particular, $M_H$ is a direct summand of the projective $R$-submodule of $M$. Hence $M_H$ is projective over $R$. In our case this implies $(M^{\otimes m})_{\Sigma_m}$ is a projective $R$-module. To see $M^{\otimes m}$ is projective over $R[\Sigma_m]$ we use the extension of scalars adjunction,
\[\hom_{R[\Sigma_m]}(M^{\otimes m}, -)\cong \hom_{R}\left(M^{\otimes m}\underset{R[\Sigma_m]}{\otimes} R, -\right)\cong \hom_{R}((M^{\otimes m})_{\Sigma_m}, -),\]
corresponding to the unique map of groups $\Sigma_m\rightarrow e$. The exactness of the above functor now follows as $(M^{\otimes m})_{\Sigma_m}$ is a projective $R$-module.
\end{proof}

\begin{proof}[Proof of Theorem~\ref{theoremfive}]
First realise $S_\ast$ by a globally flat homotopy commutative $R$-algebra $S$ with $\pi_\ast^e S\cong S_\ast$ using Theorem~\ref{theoremone}. The fact $S_0$ is a $\Z[1/n!]$-algebra implies that multiplication by $n!$ is an isomorphism on each $S_0$-module $S_q$, for all $q\in\mathbb{Z}$. We can then apply Lemma~\ref{littlehomologicalalgebra} to see $S^{\otimes m}$ is an $R$-module with $\pi_\ast^e (S^{\otimes m})$ a projective $\pi_\ast^e R[\Sigma_m]$-module. To calculate the homotopy groups of $\mathbb{L}\P_R^m S$ for every $1\leq m\leq n$ we employ Proposition~\ref{propthree}, which when evaluated at the trivial group states
\[\pi_\ast^e \mathbb{L}\P_R^m S\cong \pi_\ast^e \left(R\smp \Sigma_+^\infty E_\gl \Sigma_m\right)\underset{\pi_\ast^e R[\Sigma_m]}{\otimes} \pi_\ast^e (S^{\otimes m})\cong (S_\ast^{\otimes m})_{\Sigma_m}.\]
It follows from Lemma~\ref{littlehomologicalalgebra} again that the $R$-module $\mathbb{L}\P_R^m S$ satisfies the hypotheses of Proposition~\ref{proptwo}. Hence we obtain a natural isomorphism
\[\hom_{\ho^\gl (\mod_R)}(\mathbb{L}\P^m_R S, S)\overset{\pi_\ast^e}{\longrightarrow} \hom_{\mod_{\pi_\ast^e R}}((S_\ast^{\otimes m})_{\Sigma_m}, S_\ast).\]
Using this isomorphism we define our desired maps $h^m$, as the unique preimage of the iterated multiplication map on $S_\ast$ factored through the $\Sigma_m$-coinvariants. These maps satisfy the properties of Definition~\ref{ginfintyguys} as the iterated multiplication maps on $S_\ast^{\otimes m}$ factored through the $\Sigma_m$-coinvariants do.
\end{proof}

One can combine Theorems~\ref{theoremthree} and~\ref{theoremfive} to say that if $\eta_\ast\colon \pi_\ast^e R\rightarrow S_\ast$ is also an \'{e}tale map and $S_0$ is rational, then $S$ has a $\G_\infty$-$R$-algebra structure and an $\E_\infty$-global $R$-algebra structure. This is as close as we can get to saying $S$ has the global homotopy type of an ultra-commutative ring spectrum with the technology of this article.


\section{Examples}\label{sectionwithexamples}
Using the work above, we can show that many classical constructions in stable homotopy theory can be lifted to global homotopy theory, whilst maintaining control of the global homotopy type. We will consider localisation constructions, realisations of Galois extensions of (graded) commutative rings, some examples pertaining to periodic global complex $K$-theory, and some examples from chromatic homotopy theory.\\

For this section $R$ will denote an ultra-commutative or cofibrant $\E_\infty$-global ring spectrum (in the later case, we will use Remarks~\ref{generalisation} and~\ref{generalisation2} without reference).

\begin{example}[Localisation of algebras by an element in $\pi_\ast^e R$]\label{localisingoneelement}
Let $x\in \pi_k^e R$ be an element in the nonequivariant homotopy groups of $R$, and let $x\colon R\to \Sigma^{-k} R$ be the map representing $x$ under the representability isomorphism
\[[R, \Sigma^{-k} R]_R^\gl\cong \pi_0^e \Sigma^{-k}R\cong \pi_k^e R;\]
see Construction~\ref{constr}. Taking a fibrant replacement $\Sigma^{-k} R\to (\Sigma^{-k}R)_f$ in $\mod_R^\gl$, one can then recognise the composition of maps in $\ho^\gl(\mod_R)$
\[R\xrightarrow{x} \Sigma^{-k}R\to (\Sigma^{-k}R)_f\]
by a strict map map in $\mod_R^\gl$, which we will denote also by $x$. One can define the $R$-module $R[x^{-1}]$ as the homotopy colimit of the tower
\[R\xrightarrow{x} (\Sigma^{-k} R)_f\xrightarrow{(\Sigma^{-k}R\otimes_R x)_f}(\Sigma^{-2k} R)_f\xrightarrow{(\Sigma^{-2}R\otimes_R x)_f}\cdots,\]
where as-per-usual, $(-)_f$ denotes a fixed functorial fibrant replacement. As $R[x^{-1}]$ is defined by a filtered homotopy colimit, it is also easy to calculate $\pi_\ast^G (R[x^{-1}])$ for any $G$ in $\F$:
\[\pi_\ast^G (R[x^{-1}])\cong \colim \left(\pi_\ast^G R \xrightarrow{\cdot p_G^\ast(x)} \pi_{\ast+k}^G R\xrightarrow{\cdot p_G^\ast(x)} \pi_{\ast+2k}^G R\xrightarrow{\cdot p_G^\ast(x)}\cdots\right)\cong (\pi_\ast^G R)[p_G^\ast(x)^{-1}].\]
By inspection we see $R[x^{-1}]$ is globally flat. It is simple to place a homotopy commutative $R$-algebra structure on $R[x^{-1}]$. The unit is given by the map $R\to R[x^{-1}]$ from $R$ into the first stage of the homotopy colimit, and the multiplication map is given by the composite
\[R[x^{-1}]\smpl_R R[x^{-1}]\cong \hocolim(R\otimes_R R[x^{-1}]\xrightarrow{ x} \Sigma^{-k} R\otimes_R R[x^{-1}]\to\cdots)\]
\begin{equation}\label{multiplicationmap}\cong \hocolim(R[x^{-1}]\xrightarrow{ x} \Sigma^{-k} R[x^{-1}]\to\cdots) \xleftarrow{\cong} R[x^{-1}],\end{equation}
where the last map is the inclusion into the first stage, which is an global equivalence using the calculations above. The fact that the multiplication map (\ref{multiplicationmap}) is an isomorphism in $\ho^\gl(\mod_R)$ shows $R[x^{-1}]$ satisfies the perfect universal coefficient formula as an $R$-algebra as we shall see shortly in Lemma~\ref{disagiolemma}. The localisation part of Theorem~\ref{theoremthree} then upgrades $R[x^{-1}]$ to a globally flat $\E_\infty$-global $R$-algebra, whose global homotopy groups we totally understand. Moreover, this $\E_\infty$-global $R$-algebra structure is unique up to contractible choice.
\end{example}

Let us now prove that $R[x^{-1}]$ satisfies the hypotheses of Theorem~\ref{theoremthree}.

\begin{lemma}\label{disagiolemma}
The homotopy commutative $R$-algebra $R[x^{-1}]$ of Example~\ref{localisingoneelement} satisfies the perfect universal coefficient formula as an $R$-algebra.
\end{lemma}

\begin{proof}
Part 1 of the conditions in Definition~\ref{perfectuniversalblabla} is clear for $R[x^{-1}]$ as an $R$-algebra, as in $\ho^\gl(\mod_R)$ we have $R[x^{-1}]\smpl_R R[x^{-1}]\cong R[x^{-1}]$ as mentioned above. Using this fact, and the analogous fact in algebra, part 2 then boils down to showing the map
\[[R[x^{-1}], R[x^{-1}]]^e_R\xrightarrow{\pi_\ast^e(-)} \hom_{\mod_{\pi_\ast^e R }}((\pi_\ast^e R)[x^{-1}],(\pi_\ast^e R)[x^{-1}])\]
is an isomorphism. By two extension of scalars adjunctions, one of $\E_\infty$-global ring spectra $R\to R[x^{-1}]$ and one of graded rings $\pi_\ast^e R\to (\pi_\ast^e R)[x^{-1}]$, the above map is naturally equivalent to the isomorphism
\[\pi_0^e (R[x^{-1}])\cong[R[x^{-1}], R[x^{-1}]]^e_{R[x^{-1}]}\xrightarrow{\pi_\ast^e(-)} \hom((\pi_\ast^e R)[x^{-1}],(\pi_\ast^e R)[x^{-1}])\cong (\pi_0^e R)[x^{-1}],\]
where the hom-set above is in the category $\mod_{(\pi_\ast^e R)[x^{-1}]}$.
\end{proof}

\begin{example}[Localisation of algebras by a set in $\pi_\ast^e R$]\label{countablesetofelements}
For any countable multiplicative subset $S\subseteq \pi_\ast^e R$ one can define a globally flat $\E_\infty$-global $R$-algebra $R[S^{-1}]$ such that
\[\pi_\ast^e (R[S^{-1}])\cong (\pi_\ast^e R)[S^{-1}].\]
Indeed, one definition for such an $R[S^{-1}]$ is to enumerate $S=\{x_1, x_2,\ldots\}$, represent these elements by maps of $R$-modules as in Example~\ref{localisingoneelement}, and then define
\[R[S^{-1}]=\hocolim\left(R\xrightarrow{x_1} (\Sigma^{-k_1}R)_f\xrightarrow{(x_1x_2)} (\Sigma^{-k_2} R)_f\xrightarrow{(x_1x_2x_3)}\cdots\right),\qquad k_i=\sum_{1\leq j\leq i}|x_j|,\]
where we have suppressed some (de)suspensions of maps. The same techniques from Example~\ref{localisingoneelement} show $R[S^{-1}]$ is a globally flat homotopy commutative $R$-algebra with $\pi_\ast^e (R[S^{-1}])$ naturally isomorphic to $(\pi_\ast^e R)[S^{-1}]$, and that $R[S^{-1}]$ can be given an $\E_\infty$-global $R$-algebra structure, unique up to contractible choice. It is important here that we can commute maps representing elements in the homotopy groups of $R$ to obtain a well-defined object in $\ho^\gl(\mod_R)$.\\

Let us note that the reason we cannot extend the above result to subsets $S$ of arbitrary size is that one would like to set
\[R[S^{-1}]=\underset{\mathrm{finite}\, T\subseteq S}{\hocolim} R[T^{-1}],\]
however, using the techniques of this article, it is not clear such a filtered diagram in $\ho^\gl(\mod_R)$ can be strictified to a diagram in $\mod_R^\gl$. This is of course possible, with a more careful study of localisations, as done in \cite{me} or \cite{upcoming} in the global setting, \cite{hh} in the equivariant setting, and \cite[Section V]{ekmm} or \cite[Section 7]{ha} in the nonequivariant setting.
\end{example}

\begin{example}[Localisations of modules]\label{localisingmodules}
Given a countable multiplicative subset $S\subseteq \pi_\ast^e R$ and an $R$-module $M$, one can consider the \emph{localisation of $M$ at $S$}, defined as the global $R$-module
\[M[S^{-1}]=M\smpl_R R[S^{-1}]\simeq M\otimes_R R[S^{-1}].\]
Using the description of $R[S^{-1}]$, one obtains an alternative formula for $M[S^{-1}]$,
\[M[S^{-1}]=M\smpl_R R[S^{-1}]\simeq \hocolim \left(M\xrightarrow{\cdot x_1}(\Sigma^{-k_1}R)_f\otimes_R M \xrightarrow{\cdot(x_1x_2)}\cdots \right),\]
where we used the fact that homotopy colimits commute with derived relative smash products (up to global equivalence). One can use this formula, the global flatness of $R[S^{-1}]$, and nonequivariant flatness of localisation to obtain the calculation $\pi_\ast^G (M[S^{-1}])\cong (\pi_\ast^G M)[S^{-1}]$ for each compact Lie group $G$, where $\pi_\ast^e R$ acts on $\pi_\ast^G R$ through the algebra map induced by the unique map $p_G\colon G\to e$.
\end{example}

It is possible to generalise the above localisation examples to algebras over $\E_\infty$-global ring spectra, to localise global ring spectra at elements in \emph{equivariant} homotopy groups, and to construct localisations with ultra-commutative structure. These things are work-in-progress; see \cite{upcoming}.\\

Let us return to a counter-example from the introduction. 

\begin{example}[Global Gaussian sphere (absolute version)]\label{absolutelegaloisfailiure}
The fact that $\Z\to \Z[i]$ realises $\Z[i]$ as a projective abelian group, implies that the base change over $\pi_\ast^e \sph$ also realises $(\pi_\ast^e \sph)[i]$ as a projective $\pi_\ast^e \sph$-module. By Theorem~\ref{theoremone} we obtain a globally flat \emph{homotopy commutative ring spectrum} $\sph[i]$ realising the $\pi_\ast^e \sph$-module $(\pi_\ast^e \sph)[i]$, which is unique up to global equivalence. We claim that this homotopy commutative global ring spectrum does \textbf{not} come from an $\E_\infty$-global ring spectrum. Indeed, by the proof of Theorem~\ref{theoremone} the object $\sph[i]$ is bifibrant in $\spec^\gl$, and  by Theorem~\ref{superusefullater}, a cofibrant replacement $c\colon \sph[i]_c\to \sph[i]$ in $\spec^e$ is a global equivalence. Theorem~\ref{bifibrantendotheorem} gives us a zigzag of weak equivalences of topological operads
\begin{equation}\label{weakequivofendo}\opend_{\sph}^e(\sph[i]_c)\simeq \opend_\sph^\gl(\sph[i]).\end{equation}
For a contradiction, assume there existed a map of topological operads $\ga\colon \O\to \opend_\sph^\gl(\sph[i])$, where $\O$ is an $\E_\infty$-operad. Post-composing $\ga$ with (\ref{weakequivofendo}), and using the fact that $\O$ is cofibrant and all topological operads are fibrant, we obtain a morphism of topological operads $\O\to \opend^e_\sph(\sph[i]_c)$. This gives us an $\E_\infty$-structure on $\sph[i]_c$, which due to the projectivity of $\pi_\ast \sph[i]_c$ over $\pi_\ast \sph$ shows the natural map of $\E_\infty$-rings in $\spec^e$
\[\sph[i]_c\otimes H\Z\to H(\Z[i])\]
is an equivalence, a contradiction of \cite[Proposition 2]{svw}.
\end{example}

As in the nonequivariant case, the solution is to invert 2.

\begin{example}[Global Gaussian sphere (after inverting 2)]\label{successofinvertingtwo}
Let us now work over the $\E_\infty$-global ring spectrum $R=\sph[1/2]$ of Example~\ref{localisingoneelement}. We then have $R_\ast=(\pi_\ast\sph)[1/2]$ and using the same techniques as in Example~\ref{absolutelegaloisfailiure}, we obtain a realisation of the morphism
\[\eta_\ast\colon \pi_\ast^e \sph[1/2]\to (\pi_\ast^e \sph[1/2])[i]\]
by a globally flat $\sph[1/2]$-module spectrum $\sph[1/2, i]$. Moreover, by base change, we see the morphism $\eta_\ast$ is \'{e}tale as $\eta_0\colon\Z[1/2]\to \Z[1/2, i]$ is \'{e}tale, which is true as $\Z\to \Z[i]$ is smooth and ramified only at the prime 2. By Theorem~\ref{theoremthree}, we obtain a realisation of $S_\ast$ as a globally flat $\E_\infty$-global $R$-algebra, unique up to contractible choice, which we will denote as $\sph[1/2, i]$. Moreover, as $\sph[1/2, i]$ is globally flat over $\sph[1/2]$ then for any compact Lie group $G$ we have
\[\pi_\ast^G(\sph[1/2, i])\cong \pi_\ast^G(\sph[1/2])\otimes_{\pi_\ast^e \sph[1/2]} S_\ast\cong (\pi_\ast^G\sph)[1/2, i].\]
\end{example}

One can generalise the above example, following \cite{svw}.

\begin{example}[Adjoining roots of unity in good cases]\label{copingsvw}
Fix a prime $p$ and an integer $n\geq 1$. Suppose that $p$ is invertible inside $\pi_0^e R$ and that the $p^n$th cyclotomic polynomial 
\[\Phi_{p^n}(X)=\sum_{i=0}^{p-1}X^{ip^{n-1}}\]
is irreducible. One can then define a globally flat $\E_\infty$-global ring spectrum $R(\zeta)$ as the localisation
\[R(\zeta)=(R[C_{p^n}])\left[ \left(1-\frac{\Phi(t)}{p}\right)^{-1} \right],\]
where $R[C_{p^n}]$ is given as in Construction~\ref{grouprings}, $t$ is a generator of $C_{p^n}$, and the localisation is done \`{a} la Example~\ref{localisingoneelement}. The reason this recognises the base change over $\pi_\ast^e R$ of the map of rings $\pi_0^e R\to \pi_0^e R(\zeta)$, is due to the fact that on $\pi_0^e$, inverting the element $1-{\Phi(t)}/{p}$ is the same as taking a quotient by ${\Phi(t)}/p$, as from our hypotheses these elements are idempotents in $\pi_0^e R$; more details can be found in \cite{svw}. Furthermore, one can check that the map of graded rings $\pi_\ast^e R\to \pi_\ast^e R(\zeta)$ is \'{e}tale and realises $\pi_\ast^e R(\zeta)$ as a projective $\pi_\ast^e R$-module, so Corollary~\ref{onlyusedforglaoisextensions} states the realisation $R(\zeta)$ as an $\E_\infty$-global $R$-algebra is unique up to contractible choice. Theorem~\ref{theoremfive} states that if in addition $p!$ is invertible in $\pi_0^e R$, then $R(\zeta)$ has a $\G_p$-$R$-algebra structure as well.
\end{example}

Further generalisations of the previous two examples exist, following \cite{br}.

\begin{example}[Galois extensions of rings]\label{copingbr}
Let $G$ be a finite group and $\pi_\ast^e R\to S_\ast$ a $G$-Galois extension of graded commutative rings, so now be defined. For a finite group $G$, a \emph{$G$-Galois extension of rings} the data of a morphism of rings $A\to B$ and a $G$-action on $B$ as an $A$-algebra such that $B^G= A$ and the morphism of rings
\[\Xi\colon B\otimes_A B\to \prod_{\ga\in G} B,\qquad b_1\otimes b_2\mapsto (b_1 \ga(b_2))_{\ga},\]
is an isomorphism; see \cite[Definition 1.1.1]{br} for example. The graded case is similar. By \cite[Theorem 1.1.4]{br}, we see that $S_\ast$ is a finitely generated projective $\pi_\ast^e R$-module, so we can apply Theorem~\ref{theoremone} to obtain a globally flat homotopy commutative $R$-algebra $S$, uniquely determined in $\Ho^\gl(\mod_R)$, recognising $S_\ast$. Moreover, Theorem~\ref{theoremone} also realises the $G$-action on $S_\ast$ as a $G$-action on $S$ inside $\Ho^\gl(\mod_R)$. We note that $\pi_\ast^e R\to S_\ast$ is \'{e}tale. Indeed, if $A\to B$ is a $G$-Galois extension of rings, then using the formulation of \'{e}tale morphism as given in \cite[Definition 7.5.0.1]{ha}, we see it suffices to show the $A$-algebra multiplication map $B\otimes_A B\to B$ is the projection onto a summand. This follows though, as by definition the map $\Xi\colon B\otimes_A B\to \prod_G B$ is an isomorphism, and the mlutplication map is the composition of $\Xi$ with the projection onto the factor indexed by the identity element of $G$. Corollary~\ref{onlyusedforglaoisextensions} then shows that $S$ has an $\E_\infty$-global $R$-algebra structure, unique up to contractible choice. Moreover, this corollary also states that the natural map
\[\Map_{\calg_R^\gl}(S,S)\xrightarrow{\pi_\ast^e} \hom_{\calg_{\pi_\ast^e R}}(S_\ast, S_\ast)\]
is a weak equivalence of spaces, allowing us to lift the $G$-action on $S_\ast$ to a $G$-action on $S$ as an $\E_\infty$-global $R$-algebra. Furthermore, if $n!$ is invertible in $S_0$, then $S$ obtains a $\G_n$-$R$-algebra structure, compatible with the homotopy commutative $R$-algebra structure by Theorem~\ref{theoremfive}.
\end{example}

We now begin with two of examples involving the global complex $K$-theory spectra defined in \cite{s}.

\begin{example}[All modules over $\pi_\ast^e \mathbf{KU}$ are realisable]\label{alleasymodulesarerealisable}
In \cite[Section 6.4]{s}, Schwede constructs the ultra-commutative ring spectrum $\mathbf{KU}$, the \emph{periodic global complex $K$-theory spectrum}, which for each compact Lie group $G$ and each finite CW-complex $A$, comes with an isomorphism between the group $\mathbf{KU}_G^0(A_+)$ and the Grothendieck group of isomorphism classes of $G$-vector bundles over $A$; see \cite[Corollary 6.5.23]{s}. Moreover, the underlying nonequivariant homotopy type of $\mathbf{KU}$ is the classical complex $K$-theory spectrum; see \cite[Remark 6.4.15]{s}. This implies that $\pi_\ast^e \mathbf{KU}\cong \Z[\be^{\pm 1}]$, where $\be\in \pi_2^e \mathbf{KU}$ is the Bott element; see \cite[Construction 6.4.28]{s}. It follows that all $\pi_\ast^e\mathbf{KU}$-modules are 2-periodic, hence the data of two $\pi_0^e \mathbf{KU}$-modules, i.e., the data of two abelian groups. This implies that all graded $\pi_\ast^e \mathbf{KU}$-modules have projective dimension of 1 or less. To apply Proposition~\ref{wolbertgenerlaiseation}, we need to check a particular Tor-condition, but this follows from the facts that $\epi_0\mathbf{KU}\cong\mathbf{RU}$ and global Bott periodicity \cite[Theorem 6.4.29]{s}.\\

Indeed, as for each compact Lie group $G$, then complex representation ring $\mathbf{RU}(G)$ is a free $\Z$-module, as any finite dimensional complex $G$-represenation splits as a unique sum of simple $G$-complex representations. This shows $p_\ast^\ast\colon \pi_0^e \mathbf{KU}\to \pi_\ast^G\mathbf{KU}$ views the codomain as a free module over $\Z$. Equivariant Bott periodicity states that $\pi_\ast^G\mathbf{KU}\simeq \mathbf{RU}(G)[\be^{\pm}]$, where $\be$ is the image of the classical Bott periodicity element from $\pi_2\mathrm{KU}$. In summary, we obtain the calculation
\[\tor_k^{\pi_\ast^e\mathbf{KU}}(\pi_\ast^G \mathbf{KU}, M_\ast)\cong \tor_k^{\Z[\be^{\pm}]}(\mathbf{RU}(G)[\be^{\pm 1}], M_\ast)\cong \tor_k^{\Z}(\mathbf{RU}, M_\ast)=0,\qquad k\geq 1\]
for every $\pi_\ast^e\mathbf{KU}$-module $M_\ast$ and every compact Lie group $G$.\\

This allows us to use Proposition~\ref{wolbertgenerlaiseation}, which states that every graded $\Z[\be^{\pm 1}]$-module $M_\ast$ can be realised by a (not necessarily unique) globally flat $\mathbf{KU}$-module.
\end{example}

\begin{example}[Periodic $K$-theory from connective $K$-theory]
Writing $\mathbf{ku}^c$ for the ultra-commutative ring spectrum of \emph{global connective complex $K$-theory} (see \cite[Construction 6.4.32]{s}) and $\be\in \pi_2^e \mathbf{ku}^c$ for the Bott class (called $\nu$ in \cite[p. 648]{s}), then we claim that the $\E_\infty$-global ring spectrum $\mathbf{ku}^c[\be^{-1}]$ is globally equivalent to periodic global complex $K$-theory $\mathbf{KU}$. Indeed, there is a morphism of ultra-commutative ring spectra $\mathbf{ku}^c\to \mathbf{KU}$ (see \cite[(6.4.33)]{s}) which becomes an equivalence after localising at $\be$. This example is inherently tautological, as the definition of $\mathbf{ku}^c$ requires the definition $\mathbf{KU}$ as an input anyhow. 
\end{example}

The next example uses the well-defined global homotopy type $\mathbf{MU}$ of \cite[Section 6]{s} to lift constructions from chromatic homotopy theory to global homotopy theory. The techniques used here are essentially thoses of \cite[Section V.4]{ekmm} combined with Section~\ref{globalyflatrmodules}.

\begin{example}[Global Morava $K$-theory spectra]\label{moravaktheory}
For this example, restrict to the global family $\F=\Ab$ of \emph{abelian compact Lie groups} and fix a prime $p$. Let $\MU$ be the ultra-commutative \emph{global complex cobordism spectrum} of \cite[Example 6.1.53]{s}. It is explained in \cite[Example 6.1.53]{s} that $\MU$ has the nonequivariant homotopy type of the classical complex cobordism spectrum $\mU$ found in \cite[Example III.2.4]{bluebook}. Recall Quillen's and Lazard's theorem, which combined state that
\[\pi_\ast^e \MU_{(p)}\cong \Z_{(p)}[x_1, x_2, \ldots],\qquad |x_i|=2i\]
where we can assume that $x_{p^i-1}=v_i$, where the elements $v_i$ correspond to the Hazewinkel generators. Writing $M(k)$ for the cofibre of the multiplication by $x_k$ map on $\MU_{(p)}$,
\[\cdot x_k\colon \Sigma^{2k}\MU_{(p)}\to \MU_{(p)},\qquad k\geq 0,\]
where $x_0=p$. Define, for any $n\geq 0$, the \emph{global $n$th Morava $K$-theory spectrum} $\KK(n)$ as the $\MU_{(p)}$-module
\[\KK(n)=\left(\bigotimes_{i\neq p^n-1} M(i) \right)\smpl_{\MU_{(p)}} \MU_{(p)}[v_n^{-1}], \]
where the above smash product is relative to $\MU_{(p)}$ and derived, a countably infinite relative smash product is defined as the sequential homotopy colimit of the finite stages, and the localised $\E_\infty$-global ring spectrum $\MU_{(p)}[v_n^{-1}]$ is from Example~\ref{localisingoneelement}. Analysing the nonequivariant construction given in \cite[Section V.4]{ekmm} (or \cite[Lecture 22]{lurielecturenotes}), we see $\KK(n)$ has the nonequviariant homotopy type of classical height $n$ Morava $K$-theory as $\MU$ has the nonequivariant homotopy type as $\mU$. We claim that each $M(k)$ is globally flat over $\MU_{(p)}$. To see this, we use the fact that the morphism $p_G^\ast\colon \pi_\ast^e \MU\to\pi_\ast^G \MU$ recognises the $\pi_\ast^G \MU$ as a free $\pi_\ast^e \MU$-module for all abelian compact Lie groups $G$ -- a statement which we can transfer from that for a fixed abelian compact Lie group $G$, found in \cite[Theorem 1.3]{freenessreference}, as by \cite[Example 6.1.5.3]{s} the global spectrum $\MU$ is a model for tom Dieck's equivariant bordism for a fixed compact Lie group $G$. This means the map induced by multiplication by $x_i\in \pi_\ast^e\MU_{(p)}$,
\[\cdot p_G^\ast(x_i)\colon \pi_{\ast-2i}^G\MU_{(p)}\rightarrow \pi_{\ast}^G\MU_{(p)},\]
is injective. From this, the bottom row in the commutative diagram of $\pi_\ast^e \MU_{(p)}$-modules is exact,
\renewcommand{\M}{\mathbf{M}}
\[\begin{tikzcd}
	&	{\pi_\ast^G \M\underset{\pi_\ast^e \M}{\otimes}\pi_\ast^e \Sigma^{2i} \M}\ar[r]\ar[d, "\cong"]	&	{\pi_\ast^G \M\underset{\pi_\ast^e \M}{\otimes}\pi_\ast^e \M}\ar[r]\ar[d, "\cong"]	&	{\pi_\ast^G \M\underset{\pi_\ast^e \M}{\otimes}\pi_\ast^e M(i)}\ar[d, "{\Lambda_{M(i)}^G}"]\ar[r]	&	0	\\
{0}\ar[r]	&	{\pi_\ast^G \Sigma^k \M}\ar[r, "{\cdot p_G^\ast(x_i)}"]	&	{\pi_\ast^G \M}\ar[r]	&	{\pi_\ast^G M(i)}\ar[r]	&	0
\end{tikzcd},\]
where $\M=\MU_{(p)}$ above. By the five lemma, we see $M(i)$ is globally flat over $\MU_{(p)}$. We can compute $\epi_\ast(M(i)\smpl_{\MU_{(p)}} M(j))$ for $i\neq j$ from the short exact sequence
\[0\to \epi_\ast \Sigma^k M(i)\xrightarrow{\cdot v_j} \epi_\ast M(i)\to \epi_\ast (M(i)\smpl_{\MU_{(p)}} M(j))\to 0.\]
From this we see $M(i)\smpl_{\MU_{(p)}} M(j)$ is globally flat over $\MU_{(p)}$. By induction, each finite stage in the sequential homotopy colimit defining the infinite derived relative smash product
\[\left(\bigotimes_{i\neq n} M(i) \right)\]
is globally flat over $\MU_{(p)}$. By Proposition~\ref{naturalnessandstuff} we then see the $\MU_{(p)}$-module above is also globally flat, and Example~\ref{localisingmodules} leads us to the fact that $\KK(n)$ too, is globally flat over $\MU_{(p)}$. Nonequivariant Morava $K$-theory is most useful when considered as a ring spectrum, and either \cite[Section V.4]{ekmm} or the proof of  \cite[Lecture 22, Lemma 2]{lurielecturenotes} seemlessly work in our case too, giving $\KK(n)$ the structure of a globally flat homotopy associative $\MU_{(p)}$-algebra.
\end{example}

Using the same techniques of \cite[Section V.4]{ekmm}, one can construct globally flat $\MU_{(p)}$-modules of \emph{global Brown--Peterson spectra} $\BP$, its truncations $\BP\langle n\rangle$, \emph{global height $n$ Johnson--Wilson theory} $\mathbf{E}(n)$, and \emph{global connective height $n$ Morava $K$-theory} $\mathbf{k}(n)$. We won't mention the details of these objects here, as the only way these constructions deviate from \cite{ekmm} is by using Schwede's model for $\MU$ and the adjective globally flat, and none of the realisation results of Sections~\ref{additivesection}-~\ref{ginfinity} where used. We hope these ideas could be used in combination with the recent work of Hausmann on a Quillen's theorem for compact abelian Lie groups (see \cite{hausmann2019global}) to further the study of global chromatic homotopy theory.

\newpage
\addcontentsline{toc}{section}{References}
\bibliography{references} 
\bibliographystyle{alpha}

\end{document}